\documentclass[draft]{amsart}

\usepackage{amssymb}

\usepackage[lowtilde]{url}

\theoremstyle{plain}

\theoremstyle{definition}

\numberwithin{theorem}{section}
\numberwithin{equation}{section}

\numberwithin{enumi}{equation}

\newcommand\thmcall[1]{
  \setcounter{theorem}{\value{equation}}
  \numberwithin{equation}{theorem}    
  \numberwithin{enumi}{theorem}       
  \begin{#1}
  }

\newcommand\exitthmcall[1]{
    \end{#1}
    \numberwithin{equation}{section}  
    \numberwithin{enumi}{equation}    
    \setcounter{equation}{\value{theorem}}
  }

\newcommand\enumcall[1]{
  \stepcounter{equation}
  \begin{#1}
  }

\newcommand\exitenumcall[1]{
  \end{#1}
  }

\newcommand\bdry{\partial}
\DeclareMathOperator{\cl}{cl}
\DeclareMathOperator{\interior}{int}
\DeclareMathOperator{\LH}{LH}  

\DeclareMathOperator{\ord}{ord}
\newcommand\card{\text{card}}  

\newcommand\union{\mspace{0.5mu}\mathbin{\cup}\mspace{0.5mu}}  
\newcommand\comp{\mathbin{\sim}}  
\newcommand\with{\mspace{2mu} {\circ} \mspace{2mu}} 

\newcommand\bbz{\mathbb{Z}}
\newcommand\bbr{\mathbb{R}}
\newcommand\bbc{\mathbb{C}}
\newcommand\bbh{\mathbb{H}}
\newcommand\bbp{\mathbb{P}}

\newcommand\dx{d x}
\newcommand\dy{d y}
\newcommand\dz{d z}

\newcommand\dmu{d \mu}
\newcommand\dnu{d \nu}
\newcommand\dtheta{d \theta}

\newcommand\abs[1]{\ensuremath{\vert #1 \vert}}
\newcommand\norm[1]{\ensuremath{\Vert #1 \Vert}}
\newcommand\inprod[2]{\ensuremath{\langle #1, #2 \rangle}}

\newcommand\plusperp{\ensuremath{\oplus_\perp}}

\newcommand\orthcomp{\ensuremath{\ominus_\perp}}

\newcommand\plusalg{\ensuremath{\oplus}}

\begin{document}

\date{}

\title[Beurling's Theorem for Valuation Hilbert Modules]{Beurling's Theorem for Valuation Hilbert Modules in Several Complex Variables}
\author{Charles W. Neville}

\address{1 Chatfield Drive, Apartment 136\\
West Hartford, CT 06110\\
USA}

\email{chip.neville at gmail.com}


\thanks{Dedicated to the memory of the late Lee Rubel.  Lee, this one's for you.}

\begin{abstract}
We develop a theory of Valuation Hilbert Modules and prove a version of Beurling's theorem for these. Then we apply our version of Beurling's theorem to obtain complete descriptions of the closed invariant subspaces of a number of Hilbert spaces of analytic functions in several complex variables, including $H^2$ of the polydisk the ball, and bounded symmetric domains, and weighted $A^2$ spaces on complex analytic manifolds.
\end{abstract}

\maketitle


\setcounter{section}{0}

{\sc Part I. Introduction.}

\section{Introduction.}\label{sec1}

In a celebrated 1949 paper, Arne Beurling completely characterized closed shift invariant subspacess of $H^2$ of the unit disk.  The question of how to extend Beurling's theorem to $H^2$ of the polydisk was raised as far back as 1964 by Henry Helson \cite{Hel1}, p.\@ 25 ff.  This more than 50 year old problem was solved in a 2017 preprint and 2019 paper by Maji, Aneesh, Sankar, and Sarkar, \cite{MASS1}, \cite{MASS2}.

In this paper, we shall present an abstract version of Beurling's theorem for valuation Hilbert modules.  From this theorem, we shall deduce versions of Beurling's theorem for a number of Hilbert spaces of analytic functions in several complex variables.  These spaces include $H^2$ of the polydisk, the ball, and bounded symmetric domains, as well as weighted $A^2$ spaces on complex analytic manifolds.  ($A^2$ spaces are Hilbert spaces of analytic functions whose absolute values are square integrable with respect to weighted volume measures. They are also called Bergman spaces.)

The general theory of valuation Hilbert modules presented here, and the generalizations of  Beurling's theorem for these, $H^2$ of the ball, bounded symmetric domains, and weighted $A^2$ spaces on complex analytic manifolds are new.  As explained above, a generalization of Beurling's theorem for the polydisk is not new, but our generalization is considerably different from that by Maji, Aneesh, Sankar, and Sarkar.

\section{Preview.}\label{sec2}
In 1961, Halmos \cite{Hal1} presented a simple \emph{wandering subspace} proof of Beurling's theorem for $H^2$ of the unit disk.  Hoffman's book \cite{Hoff1} contains an excellent exposition of this material.  If $V$ is a closed, shift invariant subspace of $H^2$, a wandering subspace is a closed subspace $W$ of $V$ such that $V = W \plusperp z W \plusperp z^2 W \plusperp \cdots$.

Throughout this paper, $R$ will be an algebra over the complex numbers, and $\bbh$ will be a complex Hilbert space which is an $R$ module.  Suppose $\ord_\bbh$ and $\ord_R$ are valuations on $\bbh$ and $R$, respectively.  Then we shall call $(\bbh, \ord_\bbh)$ a \emph{valuation Hilbert module} over the valuation algebra $(R, \ord_R)$.  We have denoted the valuations by $\ord$ because, in applications involving analytic functions, $\ord(h)$ is typically the order of the zero of $h$ at a distinguished basepoint $a$.

Where it will cause no confusion, we shall denote both the algebra valuation on $R$ and the Hilbert module valuation on $\bbh$ by the same symbol, $\ord$.  We shall also make the gloss, whenever convenient, of denoting the valuation algebra $(R, \ord)$ by $R$ alone, and the valuation Hilbert module $(\bbh, \ord)$ by $\bbh$ alone.  And because this is a preview, we shall wait until later in the paper to give the general definitions of valuations in proper detail.

In place of Halmos's wandering subspace $W$, our generalization of Beurling's theorem to valuation Hilbert modules will involve the \emph{valuation homogeneous decomposition}.  Let $V$ be a closed subspace of $\bbh$.  The \emph{valuation subspace series} of $V$ is the sequence of closed subspaces $V = V_0 \supseteq V_1 \supseteq V_2 \cdots$, where $V_k = \{ h \in V \colon ord(h) \geq k \}$.  The \emph{valuation homogeneous decomposition} of $V$ is the orthogonal direct sum $V = W_0 \plusperp W_1 \plusperp W_2 \cdots$, where $W_k = V_k \orthcomp V_{k+1}$.  (Compare with Halmos's $W = V \orthcomp zV$.)  The name \emph{valuation homogeneous decomposition} arises from the fact the $\ord$ function is constant (or $0$) on each valuation homogeneous component $W_k$.

Because this is a preview, we have not offered any proofs.  Nor have we specified the necessary and sufficient conditions on the valuation homogeneous decomposition that constitute our generalization of Beurling's theorem to valuation Hilbert modules.  These will be given below in the body of our paper in {\sc Part II}.

To extend Beurling's theorem to the Hilbert spaces of analytic functions mentioned in the introduction, namely $H^2$ of the polydisk, the ball, bounded symmetric domains, and weighted $A^2$ spaces on complex analytic manifolds, we shall show that each of these, with a suitably defined $\ord$ function and complex algebra $R$, is a valuation Hilbert module over $R$.  Then we shall be able to read off a generalization of Beurling's to each of these from our generalization of Beurling's theorem for valuation Hilbert modules.

Finally, for short throughout this paper, we shall refer to our generalization of Beurling's theorem for valuation Hilbert modules as our \emph{abstract Beurling's theorem}.

\section{A Brief History of the Problem.}\label{sec3}

The first major advance in the problem of extending Beurling's theorem to $H^2$ of the polydisk was made by Ahern and Clark in 1969.  The theorem of Ahern and Clark characterized the closed submodules of $H^2$ of the polydisk of finite codimension as the closures of polynomial ideals of finite codimension \cite{AC1}.  See also the the exposition in \cite{Ru2}.

But in 1986, Agrawal, Clark, and Douglas observed, ``Almost everyone who has thought about this topic must have considered the corresponding problem for $H^2$ of the polydisk", but that little progress had been made, except for counter-examples and the important theorem of Ahern and Clark \cite{ACD1}.

After Agrawal, Clark, and Douglas wrote this, significant progress was made in related areas (cf.\@ the papers \cite{ACD1}, \cite{AB1}, \cite{Pu1}, \cite{DPSY1}, and \cite{Guo1}).  In 2007, there was a characterization of an important class of closed submodules of $H^2$ of the polydisk with deep ties to operator theory \cite{SY1}.  And as we have seen, the problem was finally definitively solved in a 2017 preprint and a 2019 paper by Maji, Aneesh, Sankar, and Sarkar \cite{MASS1}, \cite{MASS2}.

Note and Correction: Our paper is a combination, revision, and extension of our two 2021 preprints, \cite{N3} and  \cite{N4}. When we posted these, we were unaware that the problem of generalizing Beurling's theorem to the polydisk had already been solved.  We mistakenly claimed that we were the first to solve the problem.  That honor goes, of course, to Maji, Aneesh, Sankar, and Sarkar \cite{MASS1}, \cite{MASS2}.

\

{\sc Part II. Valuation Hilbert Modules and Beurling's Theorem.}
 
\section{Notation.}\label{sec4}
$\bbc$ will be the complex numbers and $\bbc^n$ complex n-space.  $\bbr$ will be the real numbers and $\bbr_+$ the positive reals, $\{r \in \bbr : r \geq 0\}$.  $\bbz$ will be the integers and $\bbz_+$ the positive integers, $\{m \in \bbz : m \geq 0\}$.

$ \bbp_n = \bbp_n [ \bbc, z]= \bbp_n [ \bbc, z_1, \dots, z_n ] $ will be the algebra of polynomials with complex coefficients in the $n$ complex variables $ z = z_1, \dots, z_n $. We shall shamelessly identify polynomials with functions, so $ \bbp_n $ can be thought of as a space of complex-valued valued polynomial functions in the variables $ z = z_1, \dots, z_n $.

\

$R$ will be a commutative complex algebra (sometimes without left or right units) and $\bbh$ a left Hilbert module over $R$.  That is $\bbh$ will be both a Hilbert space and a left $R$ module, and the multiplication maps $m_r \colon h \mapsto rh$ will be continuous for each fixed $r \in R$.

Our theory extends to cases where $R$ is a non-commutative algebra, for example a non-commutative operator algebra, but we shall not consider these here,

\

The polydisk in $\bbc^n$ of radius $0 < r < \infty$ centered at $a = (a_1, \ldots, a_n) \in \bbc^n$ is the domain
   \begin{equation*}
       D^n (a, r) = \{z = (z_1, \ldots, z_n) \in \bbc^n \colon \abs{z_j - a_j} <  r \text{ for }  j = 1, \ldots n\}.
   \end{equation*}
The unit polydisk iin $\bbc^n$ is the domain $D^n (0, 1)$.

The ball in $\bbc^n$ of radius $0 < r < \infty$ centered at $a = (a_1, \ldots, a_n) \in \bbc^n$ is the domain 
   \begin{equation*}
       B^n(a, r) = \{z = (z_1, \ldots, z_n) \in \bbc^n \colon \sum\limits_{j=1}^n \abs{z_j - a_j}^2 <  r \}.
   \end{equation*}
The unit ball in $\bbc^n$ is the domain $B^n(0, 1)$.

Orthogonal direct sums will be denoted by $\plusperp$, and orthogonal complements will be denoted by $\orthcomp$.  Algebraic direct sums will be denoted by $\plusalg$, and set theoretic complements will be denoted by $\comp$.

$\bbh$ will be a complex Hilbert space, that is the scalars will be complex numbers. Our theory applies equally well to real Hilbert spaces, but we shall not consider these here.

If $X$ is a subset of a larger topological space $Y$, we shall denote the closure of $X$ by $\cl X = \cl_Y X$, the interior of $X$ by $\interior X = \interior_Y X$, and the topological boundary of $X$ by $\bdry X = \bdry_Y X$.  Finally, $\card(X)$ will denote the cardinality of the set $X$.

\

\section{The General Theory of Valuation Hilbert Modules.}\label{sec5}

The order of the zero of an analytic function at a point $a$ in its domain is very similar to a valuation in abstract algebra.  This motivates our definition of analytic valuations, and using $\ord$ to denote them.:

\thmcall{definition}\label{def5.1}
An \emph{analytic valuation} on a complex algebra $R$ is a function $\ord _R \colon \! R \mapsto \bbz_+ \union \{\infty\}$ such that for all $r$ and $s \in R$,
\enumcall{enumerate}
  \item $\ord_R(r) = 0$ if $r$ is a left or right unit of $R$. \label{item5.1.1}
  \item $\ord_R(r) = \infty$ if and only if $r = 0$. \label{item5.1.2}
  \item $\ord_R(r s) \geq \ord_R(r) + \ord_R(s)$.  \label{item5.1.3}
  \item $\ord_R(\lambda r) = \ord_R(r)$ for $\lambda \in \bbc, \lambda \neq 0$ \label{item5.1.4}
  \item $\ord_R(r + s) \geq \min(\ord_R(r), \ord_R(s))$  \label{item5.1.5}
\exitenumcall{enumerate}
\exitthmcall{definition}
Of course, condition (\ref{item5.1.1}) is satisfied vacuously if $R$ does not have left or right units..

Here and throughout, we follow the usual conventions with respect to $\infty$: $m < \infty$ for all $m \in \bbz$, and for such $m$,  $m \cdot \infty = \infty \cdot m = \infty$ if $m \neq 0$.  Further, $\infty \cdot 0 = 0 \cdot \infty = \infty$ and $\infty \cdot \infty = \infty$.

\thmcall{definition}\label{def5.2}
A \emph{valuation algebra} is an ordered pair $(R, \ord_R)$, where $R$ is a complex algebra and $\ord_R$ is an analytic valuation on $R$.
\exitthmcall{definition}

\thmcall{definition}\label{def5.3}
Let $(R, \ord_R)$ be a valuation algebra, and let $\bbh$ be a complex Hilbert space which is a left $R$ Hilbert module.  A \emph{Hilbert module valuation} on $\bbh$ (with respect to $(R, \ord_R)$) is a function $\ord_\bbh \colon \bbh \mapsto \bbz_+ \union \{\infty\}$ such that for all $h$, $h_1$ and $h_2 \in \bbh$, and for all $r \in R$,
\enumcall{enumerate}
  \item $\ord_\bbh(h) = \infty$ if and only if $h = 0$. \label{item5.3.1}
  \item $\ord_\bbh(r h) \geq \ord_R(r) + \ord_\bbh(h)$  \label{item5.3.2}
  \item $\ord_\bbh(\lambda h) = \ord_\bbh(h)$ for $\lambda \in \bbc, \lambda \neq 0$. \label{item5.3.3}
  \item $\ord_\bbh(h_1 + h_2) \geq \min(\ord_\bbh(h_1), \ord_\bbh(h_2))$. \label{item5.3.4}
  \item The $\ord_\bbh$ function is upper semi-continuous on $\bbh$. \label{item5.3.5}
\exitenumcall{enumerate}
\exitthmcall{definition}

\thmcall{definition}\label{def5.4}
   Let $(R, \ord_R)$ be a valuation algebra. A \emph{valuation Hilbert module} over $(R, \ord_R)$ is an ordered pair $(\bbh, \ord_\bbh)$, where $\bbh$ is a left complex Hilbert module over $R$ and $\ord_\bbh$ is a Hilbert module valuation on $\bbh$ with respect to $(R, \ord_R)$.
\exitthmcall{definition}

   Where it will cause no confusion, we shall denote both the algebra valuation on $R$ and the Hilbert module valuation on $\bbh$ by the same symbol, $\ord$.  We shall also make the gloss, whenever convenient, of denoting the valuation algebra $(R, \ord)$ by $R$ alone, and the valuation Hilbert module $(\bbh, \ord)$ by $\bbh$ alone.  Finally, as indicated above, we note that the concepts of a valuation algebra and a valuation Hilbert module apply equally well to \emph{real} algebras and \emph{real} Hilbert modules. But we shall not consider these here.  Our focus will be exclusively on \emph{complex} valuation algebras and complex Hilbert modules.
   
   The properties of semicontinuous functions are given in, for example, references \cite{Wiki1} and  \cite{H1}.  The easiest way to show that a given $\ord$ function is upper semicontinuous  to use the $\limsup$ definition of upper semicontinuity, $\ord(h) \geq \limsup \ord(h_k)$.

   The simplest example of a valuation Hilbert module is $H^2$ of the unit disk:

\thmcall{example}\label{example5.5}
    Let $\bbh$ be $H^2$ of the unit disk, let $R$ be the algebra of polynomials in one  complex variable $z$, and let $\ord(h)$ be the order of the zero of $h$ at $0$ for $h \in \bbh$.  Similarly, let $ord(r)$ be the order of the zero of $r$ at $0$ for $r \in R$..  Clearly, the $\ord$ function satisfies the conditions to be a Hilbert space valuation, respectively analytic algebra valuation, except possibly for upper semicontinuity.

    To show that the $\ord$ function is upper semicontinuous, let $h_k$ be a sequence in $\bbh$ converging in $H^2$ norm to $h$. Clearly, $\ord(h) \geq \limsup \ord(h_k)$, so the $\ord$ function is, in fact, upper semicontinuous as well.  Thus $(\bbh, \ord)$ is a valuation Hilbert module over the valuation algebra $(R, \ord)$.
\exitthmcall{example}

\thmcall{example}\label{example5.6}
   Let $R$, $\bbh$, $\ord$, $h$, and $h_k$ be as in the previous example.  If $h(0) \neq 0$, $\ord( h_k )$ converges to $\ord( h ) = 0$, so $\ord$ is continuous at $h$. But if $h$ has a zero at $0$, then $\ord( h ) > 0$. Choose $h_k = h + 1/k$ and observe that $h_k$ converges to $h$, but $\ord( h_k ) = 0$ does not converge to $\ord( h )$. Thus $\ord$ is not continuous at $h$.
\exitthmcall{example}

In section \ref{sec8} below on \emph{Analytic Hilbert Modules and Standard Order Functions}, we shall discuss the order of the zero of an analytic function of several complex complex variables, as well as conditions under which $\ord(h) =$ the order of the zero of $h$ at a point $a$ in the domain of $h$ is an order function.

\section{Invariant Subspaces of Valuation Hilbert Modules.}\label{sec6}

From now on, unless the contrary is stated, $V$ will be a \emph{closed} subspace of $\bbh$.

\thmcall{definition}\label{def6.1}
The subspace $V$ of $\bbh$ is $R$ \emph{invariant} if $R V \subseteq V$.  Usually, we shall simply say an $R$ invariant subspace of $V$ is \emph{invariant}.
\exitthmcall{definition}

The subalgebra $R_1 = \{ r \in R : \ord(r) \geq 1 \}$ will be particularly important in this paper.  In fact, the natural form of our abstract Beurling's theorem will concern $R_1$ invariant subspaces.

\thmcall{definition}\label{def6.2}
Let $R$ be a valuation algebra and let $\bbh$ be a valuation Hilbert module over $R$.  Then a closed subspace $V$ of $\bbh$ is $R_1$ \emph{invariant} if $R_1 V \subseteq V$.
\exitthmcall{definition}

A natural question is, when is an $R_1$ invariant subspace $V$ of $\bbh$ a left submodule, that is when is $R V \subseteq V$?  By definition, $r V \subseteq V$ for every $r \in R_1$.  The question is, what about $r V$ when $r \in R \comp R_1$?  Conditions which guarantee that $r V \subseteq V$ are given by the following definition and proposition:

\thmcall{definition}\label{def6.3}
Let $R$ be a valuation algebra with unit element (ie.\@ identity element) $1_R$.  Then $R$ is a \emph{unital} valuation algebra if
   \begin{equation}\label{eqn6.3.1}
        R = \{ \lambda 1_R : \lambda \in \bbc \} \union R_1
    \end{equation}
\exitthmcall{definition}

\thmcall{proposition}\label{prop6.4}
Let $R$ be a valuation algebra and let $\bbh$ be a left Hilbert module over $R$.  Then every $R_1$ invariant subspace of $\bbh$ is a left $R$ submodule of $\bbh$ if either of the following conditions hold:
\enumcall{enumerate}
\item  $R = R_1$, or \label{item6.4.1}
\item  $R$ is unital. \label{item6.4.2}
\exitenumcall{enumerate}
\exitthmcall{proposition}

\begin{proof}
Clear.
\end{proof}

A simple example of a valuation algebra satisfying condition \ref{item6.4.1} is

\begin{equation*}
    R = \{ f \in \bbp_n [ \bbc, z ] : f(0) = 0 \}
\end{equation*}

Two iconic examples of valuation algebras satisfying condition \ref{item6.4.2} are the algebra of polynomials $\bbp_n[\bbc, z]$ and the algebra $H^\infty(\bbc, (M, a))$ of bounded complex-valued analytic functions on a complex manifold $M$ with basepoint $a$.  In the first case, $\ord(f) =$ the order of the zero of $f$ at $0$, and in the second $\ord(f) =$ the order of the zero of $f$ at $a$.

\section{Decompositions of Valuation Hilbert modules.}\label{sec7}

The valuation on a valuation Hilbert module allows us to decompose it and its subspaces in useful ways.  From now on, as stated earlier, $V$ will be a subspace of $\bbh$, though for a moment, it might not be closed.  To avoid trivial issues, unless otherwise stated, $V \neq \{ 0 \}$.

\thmcall{definition}\label{def7.1}
The valuation ideal series for $R$ is the decreasing sequence
\begin{equation}\label{eqn7.1.1}
  R = R_0 \supseteq R_1 \supseteq R_2 \supseteq \cdots
\end{equation}
where
\begin{equation}\label{eqn7.1.2}
  R_k = \{ r \in R : \ord(r) \geq k \}.
\end{equation}

The valuation subspace series for $V$ is the decreasing sequence
\begin{equation}\label{eqn7.1.3}
  V = V_0 \supseteq V_1 \supseteq V_2 \supseteq \cdots
\end{equation}
where
\begin{equation}\label{eqn7.1.4}
  V_k = \{ h \in V : \ord(h) \geq k \}.
\end{equation}
\exitthmcall{definition}

The names \emph{valuation ideal series} and \emph{valuation subspace series} are justified by the following proposition:

\thmcall{proposition}\label{prop7.2}
Each $R_k$ is a two sided ideal, and each $V_k$ is a subspace.
\exitthmcall{proposition}

\begin{proof}
Each set $R_k$ is a two sided ideal by items \ref{item5.1.3}, \ref{item5.1.4}, and \ref{item5.1.5}.  Each set $V_k$ is a subspace by properties \ref{item5.3.2}, \ref{item5.3.3}, and \ref{item5.3.4}.
\end{proof}

\thmcall{proposition}\label{prop7.3}
If the subspace $V$ is closed, then each subspace $V_k$ is also closed.
\exitthmcall{proposition}

\begin{proof}
Suppose $V$ is closed, and let $V_k$ be one of these subspaces. Let $h_n, n = 1, 2, 3, \ldots$ be a sequence in $V_k$ converging to $h$.  Each $h_n \in V_k$, so $\ord(h_n) \geq k$. Since $V$ is closed and $V_k \subseteq V$, $h \in V$.  Because the $\ord$ function is upper semicontinuous, $\ord(h) \geq \limsup_n\nolimits  h_n$. Further, $\ord(h_n) \geq k$ because $h_n \in V_k$, so $\ord(h) \geq k$.  Hence $h \in V_k$ by definition \ref{eqn7.1.4}, so $V_k$ is closed. 
\end{proof}

From now on, we shall reinstate the assumption that, unless otherwise noted, $V$ is a \emph{closed} subspace of $\bbh$.

\thmcall{definition}\label{def7.4}
A non-zero linear subspace $W$ of $V$ is \emph{valuation homogeneous} if the $\ord$ function for $V$ is constant on $W \comp \{0\}$.
\exitthmcall{definition}

\thmcall{definition}\label{def7.5}
The \emph{valuation homogeneous decomposition} of $V$ is the orthogonal direct sum
\begin{equation}\label{eqn7.5.1}
  W_0 \plusperp W_1 \plusperp W_2 \plusperp \, \cdots
\end{equation}
where
\begin{equation}\label{eqn7.5.2}
  W_0 = V_0 \orthcomp V_1, \; W_1 = V_1 \orthcomp V_2, \; W_2 = V_2 \orthcomp V_3, \, \ldots
\end{equation}
and $V = V_0 \supseteq V_1 \supseteq V_2 \supseteq \, \cdots$ is the valuation subspace series for $V$.  The orthogonal subspaces $W_0$, $W_1$, $W_2, \, \dots$ are called the \emph{valuation homogeneous components} of $V$.

The \emph{valuation homogeneous decomposition} of $\bbh$ itself will be denoted by
\begin{equation}\label{eqn7.5.3}
  H_0 \plusperp H_1 \plusperp H_2 \plusperp \, \cdots
\end{equation}
\exitthmcall{definition}

\thmcall{proposition}\label{prop7.6}
Let $W = W_0 \plusperp W_1 \plusperp \plusperp W_2 \, \cdots$.  Then $W$ equals all of $V$.
\exitthmcall{proposition}

\begin{proof}
By construction, $W \subseteq V$.  We shall show that if $h \in V$ and $h \perp W$, then $h = 0$.  This will prove $W = V$.

So let $h \in V$ and $h \perp W$.  Then $h \perp W_0 \plusperp W_1 \plusperp \cdots \plusperp W_m$, for $m = 0, 1, 2, \ldots$.  By construction,
   \begin{equation}\label{eqn7.7}
      V = W_0 \plusperp W_1 \plusperp \cdots \, \plusperp W_m \plusperp V_{m+1}
   \end{equation}
so $h \in V_{m+1}$ for all $m \in \bbz_+$.  From the definition of a Hilbert module, we see that $\ord(h) \geq k$ for $k  = m+1, m+2, m+3, \ldots$. Thus $\ord(h) = \infty$.  By property \ref{item5.3.1} in the definition of a Hilbert module valuation, we see that $h = 0$.
\end{proof}

The next proposition justifies calling the orthogonal direct sum \ref{eqn7.5.1} \emph{the valuatiion homogeneous decomposition} of $V$.

\thmcall{proposition}\label{prop7.8}
The valuation homogeneous components of $V$ are all closed valuation homogeneous subspaces.  In fact, $\ord(h) = m$ for all non-zero $h \in W_m$.
\exitthmcall{proposition}

\begin{proof}
$W_m$ is a closed subspace because it is the orthogonal complement of a closed subspace.  Let $h \in W_m$ and suppose that $\ord(h) \neq m$.  We shall show that $h$ must equal $0$. Now $W_m = V_m \orthcomp V_{m+1}$ from the definition of the valuation homogeneous decomposition of $V$, so $h \in V_m$, but $\ord(h) \neq m$. Since $\ord(h) \neq m$, $\ord(h)$ must be $> m$, and so $h \in V_{m+1}$.  But $W_m \perp V_{m+1}$, so $h \perp h$.  Thus $h = 0$.
\end{proof}

We also need to define the valuation homogeneous decomposition of an individual element of $\bbh$.

\thmcall{definition}\label{def7.9}
Let $h \in \bbh$.  The \emph{valuation homogeneous decomposition} of $h$ with respect to $\bbh$ is the convergent series of orthogonal terms
\begin{equation}\label{eqn7.9.1}
  h = P_0^\bbh (h) \plusperp P_1^\bbh (h) \plusperp P_2^\bbh (h) \plusperp \, \cdots.
\end{equation}
The individual terms $P_k^\bbh(h)$ are called the \emph{valuation homogeneous components} of $h$ with respect to $\bbh$..

If $h \in V$, the \emph{valuation homogeneous decomposition} of $h$ \emph{with respect} to $V$ is the convergent series of orthogonal terms
\begin{equation}\label{eqn7.9.2}
    h = P_0^V (h) \plusperp P_1^V (h) \plusperp P_2^V (h) \plusperp \, \cdots
\end{equation}
The individual terms $P_k^V$ are called the \emph{valuation homogeneous components} of $h$ \emph{with respect} to $V$

In the above, $P_k^\bbh$ and $P_k^V$, are the orthogonal projections of $\bbh$ and $V$ onto $H_k$ and $W_k$ respectively.
\exitthmcall{definition}

Note that this series converges to $h$ because the valuation homogeneous decomposition of   $V$ equals all of $V$ by proposition \ref{prop7.6}.

\section{Analytic Hilbert Modules.}\label{sec8}

In this section, we shall define the class of \emph{Analytic Hilbert Modules}, which is the concrete class of valuation Hilbert modules for which our abstract Beurling's theorem, stated below in section \ref{sec13}, holds.

\

First, we remind the reader of the standard definition of the order of the zero of an analytic function on a domain in $\bbc^n$, as given by Walter Rudin (\cite{Ru2}, section 1.1.6). We shall extend his definition to connected paracompact analytic manifolds:

\thmcall{definition}\label{def8.1}
Let $\Omega$ be a connected paracompact complex analytic manifold of complex dimension $n \geq 1$ with base point $a \in \Omega$, and let $f$ be a complex-valued analytic function on $\Omega$.  Let $V \subset \Omega$ be an open coordinate neighborhood of $a$, and let $\phi : V \longrightarrow \phi(V)\subseteq \bbc^n$ be an analytic coordinate mapping. To simplify notation, let $\psi = \phi^{-1} : \phi(V) \longrightarrow V$.  Let $h$ be a complex-valued analytic function on $V$\, and let $f = h \with \psi$.  Consider the power series expansion of $f$ around $\phi(a)$ in homogeneous terms.  Then the order of the zero of $h$ at $a$ is the minimal degree of these homogeneous terms.

Although the definition of the order of the zero of $h$ at $a$ superficially appears to depend on the coordinate map $\phi$, it is clear that it is coordinate invariant.
\exitthmcall{definition}

\thmcall{definition}\label{def8.2}
Let $\Omega$ be a connected paracompact analytic manifold of complex dimension $n \geq 1$ with a  distinguished basepoint $a \in \Omega$ as in definition \ref{def8.1}. Let $R$ be a complex algebra of complex-valued analytic functions with domain $\Omega$. Define the order of $r \in R$, $\ord_R(r) = \ord (r)$, to be the order of the zero of $r$ at $a$, as defined above. 
 
Let $\bbh$ be a complex Hilbert space of complex-valued analytic functions on $\Omega$. Define the the order of $h \in \bbh$, $\ord_\bbh(h) = \ord (h)$, to be the order of the zero of $h$ at $a$, again as defined above.

We shall call the order functions $\ord = \ord_R$ and $\ord = \ord_\bbh$ \emph{standard order functions}.
\exitthmcall{definition}

\thmcall{definition}\label{def8.3}
A valuation algebra $R$ is an \emph{analytic algebra} if its order function $\ord$ is a standard order function.  $\bbh$ is an \emph{analytic Hilbert Module} on $\Omega$ over the \emph{analytic Algebra} $R$ if $\bbh$ is a left Hilbert module over $R$ under pointwise operations, if the order function $\ord$ on $\bbh$ is a standard order function, and if there exists an open neighborhood $\Omega_a \subseteq \Omega$ of the basepoint $a$ satisfying the following:  Each sequence $h_n, n = 1, 2, 3, \dots$ in $\bbh$ converging to $h$ in $\bbh$ norm also converges to $h$ uniformly on compact subsets of $\Omega_a$. In other words, the Hilbert space topology on $\bbh$ is finer than the topology of compact convergence on $\bbh$ locally on $\Omega_a$.
\exitthmcall{definition}

Note that our analytic Hilbert modules are different from Guo's analytic Hilbert modules \cite{Guo1}.  Note also that the $\ord$ function is invariant across all local coordinate systems in $\Omega$.

\thmcall{lemma}\label{lem8.4}
Let $\bbh$ be an analytic Hilbert module on $\Omega$ over the analytic algebra $R$. Then the standard order functions $\ord_\bbh = \ord$ and $\ord_R = \ord$ on $\bbh$ and $R$ are, respectively, order functions as defined in definitions\ref{def5.1} and \ref{def5.3}:
\exitthmcall{lemma}

\begin{proof}
Clearly, the $\ord$ functions on $R$ and $\bbh$ satisfy all the properties specified in definitions \ref{def5.1} and \ref{def5.3}, except possibly for upper semicontinuity.  So all we have to do is show is that the $\ord$ function on $\bbh$ is upper semi-continuous at the base point $a$.  To this end, let $V \subseteq \Omega_a$ be an open coordinate neighborhood of the basepoint $a$. Let $\phi \colon V \rightarrow \phi(V)$ be a local coordinate mapping.  Note for clarity that $\phi(V) \subseteq \bbc^n$. Without loss of generality, $\phi(a) = 0$.

Let $0 < r$, and let $D^n(0, r)$ be an $n$-dimensional polydisk centered at $0$ such that $\psi(\cl D^n(0, r)) \subset V$ and $\psi(0) = a$.

Let $h_j , j = 1, 2, 3, \dots$ be a norm convergent sequence in $\bbh$ converging to $h \in \bbh$.  By definition \ref{def8.3}, the sequence $h_j$ also converges uniformly on compact subsets of $\Omega_a$ to $h$.  Clearly, $\psi(\cl D^n(0, r))$ is such a subset, so the sequence $h_j$ converges uniformly on $\psi(\cl D^n(0, r))$ to $h$. To simplify notation, let $f_j = h_j \with \psi$, and $f = h \with \psi$.  Then the sequence $f_j$ converges uniformly on $\cl D^n(0, r)$ to $f$.

We shall freely use multi-index notation throughout the remainder of this proof, as well as the remainder of this paper.  Recall that the degree of a monomial
\begin{equation}\label{eqn8.5}
    a_m z^m = a_{m_1} z^{m_1} \cdots a_{m_n} z^{m_n}
\end{equation}
 is
\begin{equation}\label{eqn8.6}
     \abs{m} = \abs{m_1} + \cdots + \abs{m_n},
\end{equation}

Now $f$ is an ordinary analytic function, so $\ord{(f)} $ is the degree of any of the terms of the power series expansion of $f$ of lowest degree.  Let $a_m z^m$ be one of these terms,  and let
\begin{equation}\label{eqn8.7}
    \sum\nolimits_{m_j} a_{j, m_j} z^{m_j}
\end{equation}
be the power series expansion of $f_j$ for $j = 1, 2, 3, \ldots$.  Since the sequence $f_j$ converges uniformly to $f$ on the closed polydisk $\cl D^n(0, r)$, the power series coefficients of $f_j$ converge to the power series coefficients of $f$.  Thus $\abs{a_{j, m}} > (1/2) \abs{a_m}$ for all $j >$ some integer $J$.  Now $\ord(f_j)$ is the degree of the lowest degree terms in the series \ref{eqn8.7}.  Consequently, $\ord(f_j)$ is $\leq$ the degree of $a_m z^m$, that is $\ord(f_j) \leq \ord(f)$ for all $j > J$.  Since $\limsup (\ord ( f_j ))$ is the $\inf$ of the $\sup$s, we have
\begin{equation}\label{eqn8.8}
     \ord(f) \geq \limsup\nolimits_{j > J}^\infty (\ord ( f_j )) \geq \limsup\nolimits_{j=0}^\infty (\ord( f_j ))
\end{equation}
that is $\ord( f) \geq \limsup(\ord( f_j))$. Thus the standard order function $\ord_\bbh = \ord$ is upper semicontinuous on $\bbh$. Similarly the standard order function $\ord_R = \ord$ is upper semicontinuous on $R$. Consequently, both are order functions as defined in definitions \ref{def5.1} and \ref{def5.3}.
\end{proof}

\thmcall{remark}\label{rem8.9}
From now on, we shall assume that the order functions on an analytic Hilbert module and an analytic algebra are \emph{standard order functions}.
\exitthmcall{remark}

\thmcall{remark}\label{rem8.10}
The usual proof that uniform convergence on compact subsets implies the convergence of the power series coefficients depends on the Cauchy Integral Theorem.  In our case, we would need to compute the power series coefficients of an analytic function $f$ at $0$ in the closed polydisk $\cl D^n (0, r)$:
\begin{equation}\label{eqn8.10.1}
       r^m a_m = \int_0^{2\pi} f(r e^{\imath \theta}) e^{-\imath m \theta} \dtheta/(2 \pi^n).
\end{equation}
Here, the integral is repeated $n$ times, $\theta = \theta_1 \cdots \theta_n$, and $\dtheta = \dtheta_1 \cdots \dtheta_n $.
(Rudin establishes this Cauchy formula for $a_m$ in the case where $r = 1$, \cite{Ru2}, pp.\@ 4--6.)  

Usually, the Hilbert space $E$ is required to be separable (ie.\@ contains a countable dense subset) to avoid measurability issues.  However in our case, there are no measurability issues even for non-separable Hilbert spaces because $\cl D^n (0, r)$ is contained in the larger open set $\phi (V)$, and so the integrand $f$ is continuous on $\cl D^n (0, r)$.
\exitthmcall{remark}

\thmcall{definition}\label{def8.11}
In those cases where $\bbh$ may or may not be an \emph{analytic Hilbert module} and $R$ an \emph{analytic algebra}, we shall refer to them as, respectively, a \emph{general valuation Hilbert module} and a \emph{general valuation algebra}.
\exitthmcall{definition}

If $\bbh$ is a \emph{general valuztion Hilbert module} and $R$ is a \emph{general valuation algebra}, we shall say so are explicitly. And if we need to assume, as we often shall, that $\bbh$ and $R$ are, respectively, an \emph{analytic Hilbert module} and an \emph{analytic algebra}, we shall also say so explicitly.

\section{The Minimum Value Theorem}\label{sec9}

The \emph{Minimum Value Theorem} is:

\thmcall{theorem}\label{thm9.1}
Let $R$ be a valuation algebra, let $\bbh$ be a valuation Hilbert module over $R$, and let $V$ be a closed $R_1$ invariant subspace of $\bbh$.  Suppose $h \in V$. Because the $\ord$ function on $\bbh$ takes values in the discrete set $\bbz_+$, which is bounded below by $0$, the infimum of the $\ord$ function over the valuation homogeneous components of $h$ with respect to $V$ exists and is equal to $\ord(h)$.

Consequently, $\ord(h)$ is the index of the first non-zero term in the valuation homogeneous decomposition of $h$ with respect to $V$.
\exitthmcall{theorem}

This theorem turns out to be an essential part of the proof of our abstract Beurling's Theorem. We do not know whether the Minimum Value Theorem holds for general valuation Hilbert modules, but it does for analytic Hilbert Modules, as we shall show now:

\begin{proof}
Suppose $R$ is an analytic algebra and $\bbh$ is an analytic Hilbert module over $R$. Let $h \in \bbh$, and suppose $\ord(h) = m$.. We shall use the notation in the definition of the order of the zero of an analytic function given in definition \ref{def8.1}. 

To simplify notation, let $h_k = P_k^V(h)$ for $k = 0, 1, 2 \ldots$, and $f_k = h_k \with \phi$. Consider the orthogonal direct sum
\begin{equation}\label{eqn9.2}
      h = h_0 \plusperp h_1 \plusperp h_2 \plusperp \cdots
\end{equation}
and the algebraic direct sum
\begin{equation}\label{eqn9.3}
      f = f_0 \oplus f_1 \oplus f_2 \oplus \cdots
 \end{equation}
 
That the orthogonal direct sum \ref{eqn9.2} sums to $h$ is guaranteed by proposition \ref{prop7.6} and definition \ref{def7.9}. That the algebraic direct sum \ref{eqn9.3} sums to $f$ follows from equation \ref{eqn9.2} and the fact that $f _k = h_k \with \phi$. With this simplifying notation, note that he series $\ref{eqn9.2}$ and $\ref{eqn9.3}$ represent the valuation homogeneous decomposition of $h$ with respect to $V$..

Recall that $\ord(h)$ is the minimum degree of the homogeneous terms in the power series expansion of $f = h \with \phi$ around $\phi(a)$.

Because $P_k^V(h)$ is the orthogonal projection of $h$ on $W_k$, the terms of the power series for $f_k$ are all $0$ on the valuation homogeneous components $W_r$, for $r \neq k$. Thus the homogeneous term of degree $k$ is the only non-zero term in the power series expansion of $f_k$, and so $\ord(h_k)$ is the degree of the homogeneous term of degree $k$ in this power series expansion. Thus $\ord(h_k) = k$.

Because the orthogonal sum \ref{eqn9.2} converges in norm to $h$ and $\bbh$ is an analytic Hilbert module, the analytic functions
\begin{equation}\label{eqn9.4}
   g_l = \sum_0^l h_k \text{, } l = 0, 1, 2, \ldots
\end{equation}
 converge uniformly on compact subsets of $\Omega_a$ to $h$, and so the power series coefficients of
\begin{equation}\label{eqn9.5}
   g_l \with \phi =\sum_0^l f_k \text{, } l = 0, 1, 2, \ldots
\end{equation}
converge uniformly to the power series for $f$ on the closed polydisk $\cl D^n(0, 1)$.  

Let $m$ be the index of the first non-zero term in the power series expansion of $f$ on $\cl D^n(0, 1)$.  Then $m = \ord(f) = \ord(h)$. From the uniform convergence of the power series coefficients of the sum \ref{eqn9.5} to the power series coefficients of $f$, we conclude that $\ord(h)$ is the index of the first non-zero term in the valuation homogeneous decomposition of $h$ with respect to $V$, as required.
\end{proof}

\section{$R_1$ Inner Decompositions}\label{sec10}

For $H^2$ of the unit disk, closed invariant subspaces are generated by inner functions and wandering subspaces generated by these..  The analogs of inner functions for valuation Hilbert modules are the \emph{$R_1$ inner subspaces} and  \emph{$R_1$ inner decompositions} \ref{def10.1} defined below.

\

Throughout this section, $R$ will be a general valuation algebra and $\bbh$ will be a general valuation Hilbert module over $R$.

\thmcall{definition}\label{def10.1}
A closed subspace W of $\bbh$ is \emph{$R_1$ inner} if $R_1 \cdot W \perp W$.

The valuation homogeneous decomposition of $V$ is \emph{$R_1$ inner} if, for all $k$ and $m \in \bbz_+$, for all $h \in W_k$, $g \in V$, $r \in R_1$, and $\ord(r h - g) > m$, then
\begin{equation}\label{eqn10.1.1}
  r h - g \perp W_m.
\end{equation}
\exitthmcall{definition}
 
An $R_1$ inner homogeneous decomposition must satisfy more than the orthogonality condition for an $R_1$ inner subspace; it also must satisfy the ``look ahead" condition given by equation \ref{eqn10.1.1}.  In this condition, note that there is not necessarily any relation between the integers $k$ and $m$.  All the hypotheses state is that $h$ belongs to  the homogeneous component $W_k$.  Note also that we are not asserting that the element $g \in V$ exists.  All the hypotheses state is that if an element $g$ does exist as above, then $r h - g \perp W_m$.

As we shall see later, having the valuation homogeneous decomposition of $V$ be $R_1$ inner is not sufficient for $V$ to be $R_1$ invariant.  However, it is necessary:

\thmcall{proposition}\label{prop10.2}
If $V$ is $R_1$ invariant, then the valuation homogeneous decomposition of $V$ is $R_1$ inner.
\exitthmcall{proposition}

\begin{proof}
Suppose $V$ is $R_1$ invariant.  Let $h \in W_k$, $g \in V$, and $r \in R_1$.  Then $r h \in V$, so $r h - g \in V$.  If $\ord(r h - g) > m$, then $r h - g \in V_{m+1}$, so $r h - g \perp W_m$.
\end{proof}

\section{Projections.}\label{sec11}

It would be strikingly convenient if the $R_1$ inner property were a sufficient condition for $V$ to be $R_1$ invariant, but the simplest examples show that this is not so:

\thmcall{example}\label{example11.1}
Let $R = \bbp_2 [\bbc, z_1, z_2]$, and let $\bbh = H^2(\bbc, D^2)$, complex-valued $H^2$ of the $2$-dimensional polydisk $D^2 = D^2(0, 1)$.  Let the $\ord$ functions on $R$ and $\bbh$ be standard order functions.  Clearly, $(R, \ord)$ is a valuation algebra and $(\bbh, \ord)$ is a valuation Hilbert module.  Let
   \begin{equation}\label{eqn11.1.1}
      V = W_0 \plusperp W_1 \plusperp W_2 \plusperp \plusperp W_3 \plusperp \cdots
   \end{equation}
where
   \begin{equation}\label{eqn11.1.2}
      \begin{split}
          W_0 = \{ 0 \}, \; W_1 &= \LH \{z_1, z_2\}, \; W_2 = \LH \{z_1^2, z_2^2\}, \\
                  W_3 &= \LH \{z_1^3, z_1^2 z_2, z_1 z_2^2, z_2^3\}, \; \dots,
     \end{split}
   \end{equation}
and $\LH$ denotes the linear hull.  Then $V$ is closed, the decomposition \ref{eqn11.1.1} is the $R$ homogeneous decomposition of $V$, and it is $R_1$ inner.  However, $V$ is \emph{not} $R_1$ invariant because $W_2$ lacks the needed element $z_1 z_2$.  And, in fact, the closed $R_1$ invariant subspace generated by $V$ is the closed $\bbp_2[\bbc, z_1, z_2]$ submodule generated by $z_1$ and $z_2$.
\exitthmcall{example}

We need to add an additional property to obtain sufficient conditions for a closed subspace $V$ to be $R_1$ invariant.  And we need to define two sets of projections to state this property.

Henceforth, $\bbh  =\bbh_0 \supseteq \bbh_1 \supseteq \bbh_2 \cdots$ and $\bbh =  H_0 \plusperp H_1 \plusperp H_2 \cdots$ will denote, respectively, the valuation subspace series and valuation homogeneous decomposition for the valuation Hilbert module $\bbh$.  And henceforth, just to be completely clear, $V = V_0 \supseteq V_1 \supseteq V_1 \cdots$ and $V  = W_0 \plusperp W_1 \plusperp W_2 \cdots$ will denote, respectively, the valuation subspace series and valuation homogeneous decomposition for the closed non-zero subspace $V$.

Also, until we come to the \emph{third projection lemma}, lemma \ref{lem12.6}, we shall assume that $R$ and $\bbh$ are, respectively, a general valuation algebra and a general valuation Hilbert module over $R$.

\thmcall{definition}\label{def11.2}
    \begin{equation}\label{eqn11.2.1}
        \begin{split}
  P_k^\bbh \colon \bbh &\mapsto H_k, \quad Q_{k+1}^\bbh \colon \bbh \mapsto \bbh_{k+1} \\
          P_k^V \colon V &\mapsto W_k, \quad Q_{k+1}^V \colon V \mapsto V_{k+1}
        \end{split}
    \end{equation}
will denote the orthogonal projections of $\bbh$ and $V$ onto $H_k$, $\bbh_{k+1}$, $W_k$, and $V_{k+1}$, respectively. Note that
\exitthmcall{definition}

\begin{equation}\label{eqn11.3}
    \begin{split}
        P_k^\bbh + Q_{k+1}^\bbh &= I \text{ on } \bbh_k \\
        P_k^V + Q_{k+1}^V &= I \text{ on } V_k
    \end{split}
\end{equation}
\noindent where $I$ is the identity map, and

\begin{equation}\label{eqn11.4}
    \begin{split}
  \sum_{j=0}^{k} P_j^\bbh + Q_{k+1}^\bbh &= I \text{ on } \bbh \\ 
           \sum_{j=0}^{k} P_j^V + Q_{k+1}^V &= I \text{ on } V
    \end{split}
\end{equation}

\thmcall{definition}\label{def11.5}
The valuation homogeneous decomposition of $V$ has the full projection property if for each $k$ and $m \in \bbz_+$,
\begin{equation}
  P_m^\bbh (r h - g) \in P_m^\bbh (W_m)  \label{eqn10.5.1} ,
\end{equation}
whenever $r \in R_1$, $h \in W_k$, $g \in V$, and $\ord(rh - g) \geq  m$.
\exitthmcall{definition}

\thmcall{proposition}\label{prop11.6}
If $V$ is $R_1$ invariant, then the valuation homogeneous decomposition of $V$ has the full projection property.
\exitthmcall{proposition}

\begin{proof}
Suppose $V$ is $R_1$ invariant.  Let $k$ and $m \in \bbz_+$, $r \in R_1$, $h \in W_k$, and $g \in V$.  For simplicity, denote $r h - g$ by $f$, and suppose $\ord(f) \geq m$.

Because $V$ is $R_1$ invariant, $r h \in V$ so $f = r h - g \in V$.  Now there are two cases: If $ord(f) > m$, then $P_m^\bbh (f) = 0$ because $f \in \bbh_{m+1}$ and $H_m \perp \bbh_{m+1}$.  Thus $P_m^\bbh(f)  \in P_m^\bbh (W_m)$.  On the other hand, if $\ord(f) = m$, then we can write $f = f_m + F_{m+1}$, where $f_m \in W_m$ and $F_{m+1} \in V_{m+1} \subseteq \bbh_{m+1}$.  Since $F_{m+1} \in \bbh_{m+1}$ and $H_m \perp \bbh_{m+1}$, $P_m^\bbh (F_{m+1}) = 0$.  Thus $P_m^\bbh (rh - g) = P_m^\bbh (f) = P_m^\bbh (f_m) \in P_m^\bbh (W_m)$, as required.
\end{proof}

It would also be strikingly convenient if the full projection property were a sufficient condition for $V$ to be $R_1$ invariant, but the simplest examples show that this is not so:

\thmcall{example}\label{example11.7}
Let $R = \bbp_1[\bbc, z]$, and let $\bbh = H^2(\bbc, D^1)$, complex-valued $H^2$ of the unit disk $D^1$.  Let $a \in D^1$ with $a \neq 0$, and consider the simplest possible non-constant inner function not vanishing at $0$, the simple Blaschke factor
\begin{equation}
  B(z) = \frac{a - z}{1 - \bar a z} \label{eqn11.7.1}
\end{equation}
Let
\begin{equation}\label{eqn11.7.2}
  V = \{ a_0 + a_1 z + a_2 z^2 B(z) + a_3 z^3 B(z) + \cdots \: : \: \sum \abs{a_j}^2 < \infty \} 
\end{equation}
Then $V$ is a closed subspace of $H^2(\bbc, D^1)$, the one-dimensional subspaces
\begin{equation}\label{eqn11.7.3}
\begin{split}
  W_0 = \LH \{ 1 \}, W_1 &= \LH \{ z \}, \\
                     W_2 &= \LH \{ z^2 B(z) \}, W_3 = \LH \{ z^2 B(z) \}, \dots
\end{split}
\end{equation}
are closed, $R$ homogeneous, and mutually orthogonal (because $B$ is inner).  As the reader can easily verify,
\begin{equation}\label{eqn11.7.4}
  V = W_0 \plusperp W_1 \plusperp W_2 \plusperp W_3 \plusperp \cdots
\end{equation}
is the valuation homogeneous decomposition of $V$.  Clearly, the valuation homogeneous decomposition of $V$ has the full projection property.  However, $V$ is not $R_1$ invariant (ie.\@ shift invariant).  In fact, the closed shift invariant subspace generated by $V$ is all of $H^2(\bbc, D^1)$.
\exitthmcall{example}

Variations on this example, except for the part about the decomposition \ref{eqn11.7.4} being the valuation homogeneous decomposition of $V$, are standard in the one-variable theory of $H^2$ spaces.

\section{The Projection Lemmas}\label{sec12}

The \emph{first projection lemma} is:

\thmcall{lemma}\label{lem12.1}
Suppose $R$ is a general valuation algebra and $\bbh$ is a general valuation Hilbert module over $R$. Let $r \in R_1$ and $h \in W_k$.  Let $U_m$ be the subspace $U_m = W_0 \plusperp W_1 \plusperp \cdots \plusperp W_m$.  Note that $V = U_m \plusperp V_{m+1}$.

Let $f_m = g_0 + g_1 + \cdots g_m \in U_m$, where $g_0 \in W_0$, $g_1 \in W_1$, $\dots$ $g_m \in W_m$.  Suppose further that $\ord(r h - f_m) > m$.  Then the projection of $r h$ on the subspace $U_m$ is $f_m$.  Furthermore,
   \begin{equation}\label{eqn12.1.1}
        \norm{f_m}^2 = \norm{g_0}^2 + \norm{g_1}^2 + \cdots + \norm{g_m}^2 \leq \norm{r h}^2  
   \end{equation}
\exitthmcall{lemma}

\begin{proof} For any element $w \in U_m$, $w = w_0 + w_1 + \cdots + w_m$, where $w_0 \in W_0$, $w_1 \in W_1$, $\dots$, $w_m \in W_m$.  From the definition of an $R_1$ inner decomposition, definition \ref{def10.1}, and the fact that  $\ord(r h - f_m) > m$, $r h - f_m \perp W_j$ for $j = 0, 1, \dots, m$.  Thus
\begin{equation*}
     \inprod{r h - f_m}{w} = \sum_{j=0}^m \inprod{r h - f_m}{w_j} = 0
\end{equation*}
so $r h - f^m \perp U_m$.  Since, in addition, $f_m \in U_m$, $f^m$ is the projection (up to a scalar multiple) of $r h$ on $U_m$.

As for the norm of $f_m$,
\begin{equation*}
     \norm{f_m}^2 = \norm{g_0}^2 + \norm{g_1}^2 + \cdots + \norm{g_m}^2
\end{equation*}
because the $R_1$ inner decomposition of $V$ is orthogonal.  Furthermore,
\begin{equation*}
     \norm{r h}^2 = \norm{r h - f_m}^2 + \norm{f_m}^2
\end{equation*}
because $r h - f^m \perp f_m$.  Thus inequality \ref{eqn12.1.1} holds.
\end{proof}

The \emph{second projection lemma} is:

\thmcall{lemma}\label{lem12.2}
Suppose $R$ is a general valuation algebra and $\bbh$ is a general valuation Hilbert module over $R$. Let $r \in R_1$, and $h \in W_k$.  Consider the formal series
\begin{equation}\label{eqn12.2.1}
   g_0 + g_1 + g_2 + \cdots
\end{equation}
where $g_0 \in W_0$, $g_1 \in W_1$, $g_2 \in W_2$, $\cdots$, and let
\begin{equation}\label{eqn12.2.2}
   f_m = g_0 + g_1 + \cdots + g_m.
\end{equation}
Suppose $\ord(r h - f_m) > m$ for $m = 0, 1, 2, \dots$. Then the formal series \ref{eqn12.2.1}, and thus the sequence of partial sums \ref{eqn12.2.2}, both converge in norm to an element $f$ in $V$, and $r h = f$.
\exitthmcall{lemma}

\begin{proof}
Consider the formal series of positive terms
\begin{equation}\label{eqn12.3}
  \norm{g_0}^2 + \norm{g_1}^2 + \norm{g_2}^2 + \cdots
\end{equation}
By the first projection lemma, the partial sums are
\begin{equation}\label{eqn12.4}
  \norm{f_m}^2 = \norm{g_0}^2 + \norm{g_1}^2 + \cdots + \norm{g_m}^2  \leq \norm{r h}^2  
\end{equation}
Thus the increasing sequence $\norm{f_m}^2$, $m = 0, 1, 2, \ldots$, is bounded, and so Cauchy.  Because each $f_m \in V$ and $V$ is closed and thus complete, the sequence of partial sums \ref{eqn12.2.2} and thus the formal series \ref{eqn12.2.1}, both converge to an element $f \in V$.

Now by hypothesis, $\ord(r h - f_m) > m$ for all $m$.  Also, $\lim f_m = f$ as $m \rightarrow \infty$, and the $\ord$ function is upper semicontinuous.  Thus,
\begin{equation}\label{eqn12.5}
    \ord(r h - f) \geq \limsup\nolimits_{m \rightarrow \infty} \ord(r h - f_m) > \inf \sup_{m \rightarrow \infty}\nolimits  m = \infty.
\end{equation}
\noindent Thus $\ord(r h -f) = \infty$.  By property \ref{item5.3.1} of the $\ord$ function, $r h - f = 0$.
\end{proof}

The \emph{third projection lemma} is:

\thmcall{lemma}\label{lem12.6}
Suppose $R$ is an analytic algebra and $\bbh$ is an analytic Hilbert module over $R$. Let $W$ be a valuation homogeneous subspace of $\bbh$ not equal to $\{ 0 \}$, and let $m$ be the common value of the $\ord$ function on $W \comp \{ 0 \}$.  Let $L_W$ be $P_m^\bbh$ restricted to $W,$ so that
\begin{equation}\label{eqn12.6.1}
  L_W \colon W \mapsto P_m^\bbh (W) \subseteq H_m
\end{equation}
Then $L_W$ is a bounded invertible linear transformation.
\exitthmcall{lemma}

\begin{proof}
Linearity and boundedness are trivial.  In fact, $L_W = {P_m^\bbh |}_W$ and $\norm{L_W} \leq \norm{P_m^\bbh} = 1$.  To show that $L_W$ is invertible, first note that $W \subseteq W_m$. Let $h \in W$ and $h \neq 0$. Then $m = \ord(h)$, and $P_m^\bbh (h)$ is the first non-vanishing term in the valuation homogeneous decomposition of $h$ by the minimum value theorem, theorem \ref{thm9.1}, which we may apply because $\bbh$ is an analytic Hilbert module over the analytic algebra $R$. Thus $L_W (h) = P_m^\bbh (h) \neq 0$, so $\ker(L_W) = \{ 0 \}$, that is $L_W$ is invertible.
\end{proof}

In the exceptional case where $W = \{ 0 \}$, the conclusion of the proposition still holds because $h = 0$ is the \emph{only} element of $W$.

Note that $ L_W $ is invertible on its range, but not necessarily on all of $ H_m $ because its range may be a proper subset of $ H_m $.

A simple consequence of the third projection lemma is that the projection of an element $h \in W$ onto $H_m$ determines the element $h$ itself, as does the projection onto $ \bbh_{m+1}.$

An interesting question is, when does $L_W$ have a \emph{bounded} inverse?  By the closed graph theorem, the answer is yes when $P_m^\bbh (W)$ is closed.

\section{The Abstract Beurling's Theorem}\label{sec13}

As mentioned at the close of the preview in section  \ref{sec2}, we shall refer throughout this paper to our generalization of Beurling's theorem for valuation Hilbert modules as our \emph{abstract Beurling's theorem}.  Our abstract Beurling's theorem is:

\thmcall{theorem}\label{thm13.1}
Suppose $\bbh$ is an analytic Hilbert module over the analytic algebra $R$. Let $V$ be a closed subspace of $\bbh$. Then $V$ is $R_1$ invariant if and only if the $R$ homogeneous decomposition of $V$ is $R_1$ inner and has the full projection property.
\exitthmcall{theorem}

\begin{proof}
Necessity follows from proposition \ref{prop10.2}, which states that the $R$ homogeneous decomposition of a closed $R_1$ invariant subspace is $R_1$ inner, and proposition \ref{prop11.6}, which states that a closed $R_1$ invariant subspace has the full projection property.

As for sufficiency, suppose the $R$ homogeneous decomposition of $V$ is $R_1$ inner and has the full projection property.  Let $r \in R_1$ and $h \in W_k$, where
\begin{equation}\label{eqn13.2}
  V = W_0 \plusperp W_1 \plusperp W_2 \plusperp \cdots
\end{equation}
is the $R$ homogeneous decomposition of $V$.  We shall show that $r h \in V$.

First, we shall use the full projection property, definition \ref{def11.5}, to show that there are elements
\begin{equation}\label{eqn13.3}
  g_0 \in W_0, \: g_1 \in W_1, \: g_2 \in W_2, \: \dots
\end{equation}
such that for $m = 0, 1, 2, \dots$,
\begin{equation}\label{eqn13.4}
  \ord(r h - f_m) > m
\end{equation}
\noindent Here, we have denoted the sum $g_0 + \cdots g_m$ by $f_m$.  The purpose of the functions $g_m, m=0, 1, 2, \ldots$ is to cancel out the non-zero values of the function $r h$. 

Then we shall use the second projection lemma, lemma \ref{lem12.2}, to conclude that the formal sum
\begin{equation}\label{eqn13.5}
  g_0 + g_1 + g_2 + \cdots
\end{equation}
converges to an element $f \in V$, and that $r h = f$.  This will prove that $r h \in V$.

\

First, let us construct the sequence \ref{eqn13.3} by induction.  Because $r \in R_1$, $\ord(r h) > 0$.  Thus $P_0^\bbh (r h) = 0$, so we can (and must) choose $g_0 = 0$.  Now suppose $g_0, \dots, g_{m-1}$ have been constructed inductively, so that $g_j \in W_j$ and $\ord(r h - f_j) > j$ for $j = 0, \dots, m - 1$.  We shall use the full projection property of the $R$ homogeneous decomposition \ref{eqn13.2} of $V$ to construct $g_m$.

Note that $\ord(r h - f_{m-1}) \geq m$.  By the full projection property, $P_m^\bbh (r h - f_{m-1}) \in P_m^\bbh (W_m)$.  If $W_m \neq \{ 0 \}$, we may use the third projection lemma, lemma \ref{lem12.6}, to conclude there is a unique element $g_m \in W_m$ such that $P_m^\bbh (g_m) = P_m^\bbh (r h - f_{m-1})$.  In the exceptional case where $W_m = \{ 0 \}$, we can (and must) choose $g_m = 0$.  Let $f_m = f_{m-1} + g_m$.  Then $f_m = g_1 + \cdots + g_m$, and $g_m \in W_m$ by construction.

But we still must prove that $\ord(r h - f_m) > m$.  To that end, recall from the minimum value theorem, theorem \ref{thm9.1}, that $\ord(r h - f_m)$ is the index of the first non-vanishing term in the $R$ homogeneous decomposition of $r h - f_m$, or $\infty$ if $r h - f_m = 0$.  Since $\ord(r h - f_{m-1}) > m - 1$, all the terms in the $R$ homogeneous decomposition of $r h - f_{m-1}$ vanish up to and including the index $m - 1$, that is $P_j^\bbh (r h - f_{m-1}) = 0$ for $j = 0, \dots, m - 1$.  Also, because $g_m \in W_m$, $\ord(g_m) = m$ (or $\infty$ if $g_m = 0$).  Thus $P_j^\bbh (g_m) = 0$ as well for $j = 0, \dots m-1$.  Furthermore, $P_m^\bbh (r h - f_{m-1}) = P_m^\bbh (g_m)$ by construction.  Thus
\begin{equation}\label{eqn13.6}
   P_j^\bbh (r h - f_m) =  P_j^\bbh (r h - f_{m-1}) - P_j^\bbh (g_m) =0
\end{equation}
for $j = 0, \dots, m$, and so $\ord(r h - f_m) > m$.

This concludes the first part of the proof that $r h \in V$.  The second part is immediate.  Because the valuation homogeneous decomposition of $V$ is $R_1$ inner, we may apply the second projection lemma.  By this lemma, the formal series \ref{eqn13.6} converges in norm to an element $f \in V$, and $r h = f$.  Thus $r h \in V$.

Now let $h \in V$, not necessarily an element of $W_k$.  Let $r \in R_1$ as before, and let
\begin{equation}\label{eqn13.7}
  h = h_0 \plusperp h_1 \plusperp h_2 \plusperp \cdots
\end{equation}
be the $R$ homogeneous decomposition of $h$ with respect to $V$.  By the definition of the $R$ homogeneous decomposition of $h$ with resoect to $V$, definition \ref{def7.9}, each $h_k \in W_k$, and by proposition \ref{prop7.6}, the series (\ref{eqn13.5}) converges in norm to $h$.  By the definition of a Hilbert module, the linear operators $m_r \colon f \to r f$, $f \in \bbh$ are bounded in operator norm, so the series
\begin{equation}\label{eqn13.8}
  r h = r h_0 + r h_1 + r h_2 + \cdots
\end{equation}
converges in norm to $r h$.  Further, each term $r h_k \in V$, so, since $V$ is closed, $r h \in V$ as required.  Thus $V$ is $R_1$ invariant.
\end{proof}

The last part of this proof is useful enough to isolate as a separate corollary:

\thmcall{corollary}\label{cor13.9}
Let $\bbh$ be an analytic Hilbert module over the analytic algebra $R$. A closed subspace $V$ is $R_1$ invariant if and only if, for each $k \in \bbz_+$ and $r \in R_1$, $r W_k \subseteq V$.
\exitthmcall{corollary}

\begin{proof}
See the above.
\end{proof}

A proper generalization of Beurling's theorem should characterize, not just closed $R_1$ invariant subspaces of $\bbh$, but closed $R$ submodules as well.  This we can do in two important cases.  Recall from definition \ref{def6.3} that a valuation algebra $R$ is \emph{unital} if $R$ has a unit element $1_R$ and
\begin{equation}\label{eqn13.10}
  R = \{ \lambda 1_R : \lambda \in \bbc \} + R_1
\end{equation}
From proposition \ref{prop6.4} we immediately have

\thmcall{corollary}\label{cor13.11}
Let $R$ be an analytic algebra satisfying either of the following conditions:
\enumcall{enumerate}
  \item $R = R_1$, or
  \item $R$ is unital
\exitenumcall{enumerate}
Then a closed subspace $V$ of an analytic Hilbert module $\bbh$ over $R$ is an $R$ submodule of $\bbh$ if and only if the $R$ homogeneous decomposition of $V$ is $R_!$ inner and has the full projection property.
\exitthmcall{corollary}

\begin{proof}
Apply proposition \ref{prop6.4} and our abstract Beurling's theorem.
\end{proof}

\

{\sc Part III. Applications of the Abstract Beurling's Theorem.}

\section{Invariant Subspaces of $H^2$ Spaces.}\label{sec14}

In this section, we shall use our abstract Beurling's theorem to characterized the closed invariant subspaces of $H^2$ of the polydisk, the ball, and bounded symmetric domains.  Our strategy will be to prove that $H^2$ of the polydisk, the ball, and bounded symmetric domains are analytic Hilbert modules.  Thus, we may apply our abstract Beurling's theorem \ref{thm13.1}, and its corollary \ref{cor13.11}, to completely characterize their closed invariant subspaces.

Though the Hilbert spaces of analytic functions we shall consider in this paper are all complex-valued, our results could be extended to Hilbert space valued analytic functions. Therefore, we shall include the Hilbert space of values in our notation, like this for complex-valued $H^2$ of the polydisk: $H^2(\bbc, D^n)$.
 
 \
 
First consider $H^2$ of the polydisk.  Let $D^n$ be the $n$-dimensional polydisk $D^n = D^n (a, r)$, where $a \in \bbc^n$ and $0 < r  < \infty$.  Recall that an $n$-harmonic function on $D^n$ is a continuous function which is harmonic in each variable separately, and an $n$-subharmonic function is an upper semicontinuous function on $D^n$ which is subharmonic in each variable separately (cf.\@ Rudin \cite{Ru2}, pages \@ 16 and 39).  The $n$-harmonic majorant definition of $H^2$ of the polydisk is due to Walter Rudin and is contained in chapter 3 of \cite{Ru2}, and for domains in one variable was given by Walter Rudin in \cite{Ru1}.

\thmcall{definition}\label{def14.1}
$H^2$ of the polydisk $D^n = D^n(a, r)$ is the Hilbert space $\bbh = H^2(\bbc, D^n)$ of complex-valued analytic functions $h$ on $D^n$, such that the $n$-subhar\\monic function
\begin{equation}\label{eqn14.1.1}
     z \mapsto \norm{h(z)}^2, \,\, z \in D^n
\end{equation}
admits an $n$-harmonic majorant.  If $ h \in H^2(\bbc, D^n) $, let $u$ be the least $n$-harmonic majorant of the $n$-harmonic function \ref{eqn14.1.1}.  The norm of $h$ is
\begin{equation}\label{eqn14.1.2}
     \norm{h}^2 = u(a).
\end{equation}
\exitthmcall{definition}
 
Throughout this section, $R$ will either be the analytic algebra of complex-valued polynomials $ \bbp_n [\bbc, z_1, \ldots z_n] $, or the analytic algebra of complex-valued bounded analytic functions $ H^\infty ( \bbc, D^n ) $ on $ D^n $.  Note that both of these are unital so, by corollary \ref{cor13.11}, $R_1$ invariant subspaces and $R$ invariant subspaces are the same.
 
 \thmcall{theorem}\label{thm14.2}
$ \bbh $ is an analytic Hilbert module over $ R $. Therefore, our abstract Beurling's theorem \ref{thm13.1} and its corollary \ref{cor13.11} completely describe the closed invariant subspaces of $ \bbh $.
\exitthmcall{theorem}

 \thmcall{proof}
This will follow quickly from the Harnack inequalities for $n$-harmonic functions. The Harnack inequalities themselves follow immediately from the standard formula for the Poisson kernel in each variable separately.

Let $h_j,  j = 1, 2, 3, \dots$ be a norm convergent sequence in $\bbh$ converging to $h \in \bbh$. Let $f_j = h - h_j$. Then  $\norm{f_j} \rightarrow 0$ as $j \rightarrow \infty$.
Let $u_j$ be the least $n$-harmonic majorant of the $n$-subharmonic function $z \rightarrow \norm{f_j(z)}^2$.  Let $K$ be a compact subset of $D^n$.  By the Harnack inequalities, there is a constant $0 < C < \infty$ such that $u_j(w) < C u_j(a)$ for all $w \in K$.  But $\norm{f_j(w)}^2 \leq u_j (w) \leq C u_j (a) = C \norm{f_j}^2$, so $f_j$ converges uniformly on the compact set $K$ to $0$. Thus $h_j$ converges uniformly to $h$ on compact subsets of $D^n$, so the the norm topology on $\bbh$ is finer than the topology of compact convergence on $D^n$.  Further, it is clear that $\bbh$ is an analytic Hilbert module over $R$, except that we still need to show that $\bbh$ is complete, so that $\bbh$ is a Hilberrt space:

To this end, consider a Cauchy sequence $h_j$ in $\bbh$.  Then, by essentially the same argument as above, the sequence is Cauchy in the topology of uniform convergence on $D^n$, so it converges in this topology to an analytic function $h$.  Further, the norm of $h$ is finite so $h \in \bbh$, and the sequence $h_j$ converges to $h$ in norm.  Thus $\bbh$ is complete.
\end{proof}

Next, consider $H^2$ of the \emph{ball}.  Let $B^n$ be the $n$-dimensional ball $B^n (a, r)$, where $a \in \bbc^n$ and $0 < r  < \infty$.

\thmcall{definition}\label{def14.3}
$H^2$ of the ball $B^n$ is the Hilbert space $\bbh = H^2(\bbc, B^n)$ of analytic functions $h$ on $B^n$, such that the subharmonic function
\begin{equation}\label{eqn14.3.1}
     z \mapsto \norm{h(z)}^2, \,\, z \in B^n
\end{equation}
admits a harmonic majorant.  If $h \in H^2(\bbc, B^n)$, let $u$ be the least harmonic majorant of the harmonic function \ref{eqn14.3.1}.  The norm of $h$ is
\begin{equation}\label{eqn14.3.2}
     \norm{h}^2 = u(a).
\end{equation}
\exitthmcall{definition}

Note that Rudin defines $H^2(\bbc, B^n)$ as above, but in terms of $\mathcal{M}$-harmonic functions instead of ordinary harmonic functions (cf.\@ Rudin \cite{Ru4}, pp.\@ 83--85).  Here, $\mathcal{M}$ is the conformal automorphism group of the ball.  He observes that ordinary harmonic functions will do as well, but at the expense of not having results be $\mathcal{M}$ invariant.

 \thmcall{theorem}\label{thm14.4}
$\bbh$ is a Hilbert space and an analytic Hilbert module over $R$. Therefore our abstract Beurling's theorem, theorem \ref{thm13.1}, and its corollaries describe the closed invariant subspaces of $\bbh$.
\exitthmcall{theorem}

 \thmcall{proof}
The proof is almost exactly like the proof for $H^2$ of the polydisk, except that the Harnack inequalities for ordinary harmonic functions replace the Harnack inequalities for $n$ harmonic functions.  We leave the details to the reader.
\end{proof}

Now consider the case of of a product of balls.  By analogy with the case of the polydisk, we shall call such a domain a \emph{polyball}: 

\thmcall{definition}\label{def14.5}
A \emph{polyball} is a Cartesian product $B^n$ of $k$ balls of possibly different dimensions, centers, and radii:
\begin{equation}\label{eqn14.5.1}
   B^n = B^n (a, r) = B^{n_1} (a_1, r_1) \times \cdots \times B^{n_k} (a_k, r_k)
\end{equation}
where $n = (n_1, \ldots n_k)$, $a = (a_1, \ldots a_k)$, and $(r = r_1, \ldots r_k)$.
\exitthmcall{definition}

Note that the dimension $n$ is not just a simple integer, but a ``multi-integer".

\

We shall define $H^2(\bbc, B^n)$ in terms $n$-harmonic functions on the polyball $B^n$, by analogy with the definition of $H^2(\bbc, D^n)$ in terms of $n$-harmonic functions on the polydisk $D^n$\!, except that now $n$ is a multi-integer as in equation \ref{eqn14.5.1}.

\thmcall{definition}\label{def14.6}
An $n$-harmonic function on $B^n$ is a continuous function on $B^n$ which is harmonic in each variable
\begin{equation}\label{def14.6.1}
     z_{n_j} \in B^{n_j}, \,\, j = 1, \ldots, k
\end{equation}
separately, and an $n$-subharmonic function on $B^n$ is an upper-semicontinuous function on $B^n$ which is subharmonic in each variable $z_{n_j} \in B^{n_j} $ separately.
\exitthmcall{definition}

Of course, the variables $z_{n_j}$ themselves range over $n_j$ tuples of complex numbers $\in B^{n_j} \subset \bbc^{n_j}$.

\thmcall{definition}\label{def14.7}
$H^2$ of the polyball $B^n$ is the Hilbert space $\bbh = H^2(\bbc, B^n)$ of analytic functions $h$ on $B^n$, such that the $n$-subharmonic function
\begin{equation}\label{eqn14.7.1}
     z \mapsto \norm{h(z)}^2, \,\, z \in B^n
\end{equation}
admits an $n$-harmonic majorant.  If $h \in H^2(\bbc, B^n)$, let $u$ be the least $n$-harmonic majorant of the $n$-harmonic function \ref{eqn14.7.1}.  The norm of $h$ is
\begin{equation}\label{eqn14.7.2}
     \norm{h}^2 = u(a).
\end{equation}
\exitthmcall{definition}

\thmcall{theorem}\label{thm14.8}
$\bbh$ is a Hilbert space and an analytic Hilbert module over $R$. Therefore, our abstract Beurling's theorem, theorem \ref{thm13.1}, and its corollariy, corollary \ref{cor13.11}, completely describe the closed invariant subspaces of $\bbh$.
\exitthmcall{theorem}

 \thmcall{proof}
The proof using the appropriate Harnack inequalities is almost exactly like the proof for $H^2$ of the polydisk and the ball .  We leave the details to the reader.
\end{proof}

Note that the unit disk is the unit ball of dimension 1.  Thus theorem \ref{thm14.8} characterizing the closed invariant subspaces of $H^2$ of the polyball subsumes theorem \ref{thm14.2} characterizing the closed invariant subspaces of $H^2$ of the polydisk.

\

Finally, we shall consider the case of bounded symmetric domains.  Let $\Omega$ be a bounded symmetric domain in $\bbc^n$.  Recall that $\Omega$ is \emph{irreducible} if it is not the product of bounded symmetric domains of lower dimensions.  Oherwise, $\Omega$ is \emph{reducible}.  If $\Omega$ is irreducible, the definitions, theorems, and proofs for $H^2$ of $\Omega$, $\bbh = H^2(\bbc, \Omega)$, are almost exactly like the corresponding ones for $H^2(\bbc, B^n)$.  On the other hand, if $\Omega$ is reducible, these definitions, theorems, and proofs are almost exactly like the ones for the polydisk and the polyball.  In particular,the appropriate Harnack inequalities hold.  Thus we have,

\thmcall{theorem}\label{thm14.9}
$\bbh = H^2(\bbc, \Omega)$ is a Hilbert space and an analytic Hilbert module over $R$. Therefore our abstract Beurling's theorem, theorem \ref{thm13.1}, and its corollary, corollary \ref{cor13.11}, our abstract Beurling's theorem, and its corollariy, corollary completely describe the closed invariant subspaces of $\bbh$.
\exitthmcall{theorem}

\begin{proof}
See the discussion above.
\end{proof}

An extremely general example of an analytic Hilbert module would be the case of $\bbh = H^2(\bbc, \Omega, p)$, where $\Omega$ is a connected complex analytic manifold, the basepoint $a \in \Omega$, , and $R = H^\infty(\bbc, \Omega)$.

Unfortunately, there are a number of difficulties in defining $H^2(\bbc, \Omega)$.  For example, we need $\Omega$ to be a K\"{a}hler manifold (or compact) in order to guarantee that analytic functions are harmonic.  Also, the Laplacian (the Laplace Beltrami operator) often has variable coefficients, so we cannot simply use the well known potential theory for domains in $\bbc^n$. 

\

There is, however, one circumstance where the difficulties mentioned above are already worked out, and that is the the case where the the (complex) dimension of $\Omega$ is $1$, in other words where $\Omega$ is a Riemann surface with basepoint $a$..  In fact, this is the case where Walter Rudin first introduced the harmonic majorant definition of (in our notation) $H^p (\bbc, \Omega)$ for $1 \leq p < \infty$ \cite{Ru1}.  While the ordinary Laplacian is not invariant under coordinate changes, ordinary functions harmonic in local coordinates are.  Furthermore, analytic functions are harmonic.  Thus, if we define $H^2(\bbc, \Omega)$ as the space of analytic functions $h$ on $\Omega$ such that the subharmonic function
 \begin{equation*}
     z \mapsto  \norm{h(z)}^2, \,\, z \in \Omega
 \end{equation*}
 has a harmonic majorant, and if we norm $H^2(\bbc, \Omega)$ by 
 \begin{equation*}
     \norm{h}^2 = u(a),
\end{equation*}
where $u$ is the least harmonic majorant of this subharmonic function, we have theorem \ref{thm14.10}:

\thmcall{theorem}\label{thm14.10}
Let $\Omega$ be a Riemann surface with basepoint  $a \in \Omega$, $\bbh = H^2(\bbc, \Omega)$, and $R = H^\infty(\bbc, \Omega)$,  Then $\bbh$ is an analytic Hilbert module over the analytic algebra $R$. Therefore, our abstract Beurling's theorem, theorem \ref{thm13.1}, and its corollariy, corollary \ref{cor13.11} completely describe the closed invariant subspaces of $\bbh$.
\exitthmcall{theorem}

For infinitely connected Riemann surfaces $\Omega$ satisfying certain geometric conditions, the author in \cite{N1} and \cite{N2} presented a different description of the closed submodules of $H^p(\bbc, \Omega)$ for $1 \leq p < \infty$, and the $\beta$ closed ideals of $H^\infty (\bbc, \Omega)$.  This description used multiple valued inner functions $I$ with a single valued absolute value $\abs{I}$.  Later, Hasumi, \cite{Has1} and \cite{Has2}, using the language of analytic vector bundles instead of multiple valued inner functions, but not including or extending the work on the geometric conditions on $\Omega$, proved similar theorems for $H^p$ and $L^p$ spaces.  Theorem \ref{thm14.10} generalizes all these results for $H^2$ spaces.

\section{Invariant Subspaces of $A^2$ Spaces.}\label{sec15}

Let $\Omega$ be a domain in $\bbc^n$, and let $\dmu$ be (real) $2 n$-dimensional volume measure on $\bbc^n$. The space $A^2(\bbc, \Omega, \dmu)$ is the space of all complex-valued analytic functions $h$ on $\Omega$ such that
\begin{equation*}
    \int_\Omega \norm{h(z)}^2 \dmu (z) < \infty .
\end{equation*}
We shall refer to these spaces interchangeably as $A^2$ spaces and \emph{Bergman} spaces.

\

A celebrated result of Aleman, Richter, and Sundberg from the 1990s established a Beurling type theorem characterizing the closed $ \bbp_1 [ \bbc, z ]  $ submodules of $ A^2 (\bbc, D^1(0, 1 ), \dmu ) $ \cite{ARS1}, where $ \dmu $ is real $2$-dimensional area measure.  Their proof was quite difficult, and famously involved biharmonic functions.  Since then, considerably simpler proofs have were given by Shimorin (\cite{Sh1}, \cite{Sh2}.  Shimorin has extended the ARS theorem to weighted $ A^2 $ spaces, $ A^2 ( \bbc, D^1(0, 1), k\dmu ) $ of the unit disk, albeit with restrictions on the weights $ k $ such as circular symmetry and a subharmonicity condition (\cite{Sh1}, \cite{Sh2}).  And recently, Ball and Bolotnikov gave an elegant and surprising characterization of the closed $ \bbp_n [ \bbc, z ] $ submodules of weighted $A^2$ spaces of the unit disk with the same sort of restrictions on the weight functions as in Shimorin's work. The Ball-Bolotnikov papers are particularly notable in that their characterization involves analytic inner functions rather than biharmonic functions, and has important connections to operator theory (\cite{BaBo1} and especially \cite{BaBo2}).

\

In this section, we shall generalize the ARS and Shimorin theorems to several complex variables and to connected paracompact analytic manifolds.

\

Throughout this section, $\Omega$ will be a connected, paracompact analytic manifold with basepoint $a$.  Let $\dmu$ be a (real) $2 n$-dimensional volume measure on $\Omega$.  Let $k$ be an admissible weight function on $\Omega$.  (We shall discuss the existence of the measure $\dmu$ and give the definition of an admissible weight function momentarily.)  Let $R = H^\infty(\bbc, \Omega)$,  the algebra of bounded complex-valued analytic functions on $\Omega$.  In the case where $\Omega \subseteq \bbc^n$, let $R$ be either $\bbp_n[\bbc, z_1, \dots, z_n]$ or $H^\infty(\bbc, \Omega)$.

\thmcall{definition}\label{def15.1}
The weighted Bergman space on $\Omega$, $\bbh = A^2(\bbc, \Omega, k\dmu)$, is the space of all complex-valued analytic functions $h$ on $\Omega$ such that the weighted $A^2$ norm
  \begin{equation}\label{eqn15.1.1}
     \norm{h}^2 = \int_\Omega \norm{h(z)}^2 k(z) \dmu (z) < \infty .
  \end{equation}
\exitthmcall{definition}
Define the order functions $\ord(r)$ and $\ord(h)$ to be, as usual, standard order functions.

\thmcall{theorem}\label{thm15.2}
    $\bbh = A^2(\bbc, \Omega, k \dmu)$ is a complex-valued Hilbert space and an analytic Hilbert module over $R$. Therefore our abstract Beurling's theorem,  theorem \ref{thm13.1},  and its corollary, corollary \ref{cor13.11}, completely describe the closed invariant subspaces of $\bbh$.
\exitthmcall{theorem}

In the special case of unweighted Bergman spaces defined on domains in $\bbc^n$, theorem \ref{thm15.2} follows immediately from the following lemma:

\thmcall{lemma}\label{lem15.3}
  Let $V$ be a domain in $\bbc^n$.  Let $\bbh = A^2(\bbc, V, \dmu))$, and let $K$ be a compact subset of $V$. Then there is a constant $0 < C_K < \infty$ depending on $V$ and $K$, such that
   \begin{equation}\label{eqn15.3.1}
     \sup_{z \in K} \abs{h(z)} \leq C_K \norm{h} \text{   all  }  h \in A^2 (\bbc, V, \dmu) \text{  (\cite{Kr1}, lemma 1.4.1).}
   \end{equation}
Thus a sequence $h_k, k \in \bbz_+$ converging in norm to $h$ in $\bbh$ also converges to $h$ in the topology of compact convergence on compact subsets of $V$.  Hence $\bbh$ is  an analytic Hilbert module over the analytic algebra $R$.
\exitthmcall{lemma}  

\

To return to volume measures and admissible weight functions, recall that a real $2n$-dimensional volume measure on a connected paracompact complex analytic manifold $\Omega$ is a measure which, in local coordinates $z = z_1, \dots, z_n$ on a typical open coordinate neighborhood $V$, is $\dmu(z) = c(z) \dx_1 \dy_1 \cdots \dx_n \dy_n $.  Here, $c$ is Lebesgue measurable (with respect to the measure $\dx_1 \dy_1 \cdots \dx_n \dy_n$ on $V$), and bounded {a.}{e.}~on $V$. Further, for some open coordinate neighborhood $V$ of the basepoint $p$, $c >$ some $\delta > 0$ on $V$ {a.}{e.}   A suitable volume measure always exists on $\Omega$, and in fact there are many of them.  For example, if $G$ is a Riemannian metric on $\Omega$, the Riemannian volume element $\dmu = \sqrt{ \text{det}(G)} \, \dx_1 \dy_1  \cdots  \dx_n \dy_n$ is such a measure \cite{Gr1}.

\thmcall{definition}\label{def15.4}
A real valued Lebesque measurable function $k$ is an admissible weight function on $\Omega$ if it is a $2n$-dimensional volume measure on $\Omega$, and if it is locally bounded away from $0$ a.e.\@ on $\Omega$.  To be more specific, $k$ is locally bounded away from $0$ a.e.\@ if, for each $a \in \Omega$,  there exists a connected open coordinate neighborhood $V$ of $a$, such that $k$ is bounded away from $0$ a.e.\@ on $V$,  that is if there exists a constant $\epsilon_V >  0$ such that $k > \epsilon_V$ a.e.\@ on $V$.
\exitthmcall{definition}

Note that our admissible weight functions are considerably more general than Shimorin's weight functions.

\begin{proof}[Proof of Theorem \ref{thm15.2}]
Let $V$ be an open coordinate neighborhood as in definition \ref{def15.4}, and $\phi \colon V \rightarrow U$ be an analytic chart. Note for clarity that $U$ is a domain in $\bbc^n$. Without loss of generality, $\phi(p) = 0$. To simplify notation, let $\psi = \phi^{-1} \colon U \rightarrow V$.

In the notation of definition \ref{def15.4}, let $0 < \epsilon \leq k(z)$ a.e.\@ on $V$, and let $K$ be a compact subset of $V$. Let $C_{\phi(K)}$ be as in lemma \ref{lem15.3}, with $K$ replaced by the compact subset $\phi(K)$ in $U$. Let $C_K = C_{\phi(K)}$.  The theorem now follows from the chain of inequalities,
\thmcall{equation*}
  \begin{split}
\epsilon \sup_{z \in K} \norm{h(z)}^2_E &\leq \epsilon \,C_K\int_U \norm{h \with \psi(z)}_E^2 \dnu(z) = \epsilon \,C_{K} \int_V \norm{h(z)}_E^2 k(z) \dmu(z) \\
          &\leq \epsilon \,C_K \int_\Omega \norm{h(z)}_E^2 k(z) \dmu(z) = \epsilon \,C_K \norm{h}_E^2 .
   \end{split}
\exitthmcall{equation*}

Here, the left hand inequality follows from lemma \ref{lem15.3}, and the passage from the integral over $V$ to the integral over $\Omega$ follows from the monotoneity of the Lebesgue integral for positive functions.

Thus the Hilbert space topology on $\bbh$ is finer than the topology of compact convergence, and so $\bbh$ is an analytic Hilbert module over the analytic algebra $R$.  That $\bbh$ is actually a Hilbert module over $R$ also follows quickly from the fact that the topology of compact convergence is finer than the norm topology.
\end{proof}

We shall close this section with a somewhat surprising observation:  Our theory of valuation Hilbert spaces, and our theorem \ref{thm15.2}, together with our abstract Beurling's theorem and its corollary, give new results on invariant subspaces, as well as sometimes simpler proofs, even in the classical case of one complex variable.

\section{Paradox lost.}\label{sec16}

The title refers, not to Milton's great poem, but to an apparent paradox which arises from the Aleman, Richter, and Sundberg (ARS) characterization of the closed invariant subspaces of the (unweighted) Bergman space of the unit disk in $1$ complex variable, that is $A^2( \bbc, D^1(0, 1), \dz \dy)$ in our notation.

\thmcall{theorem}[Aleman, Richter, and Sunberg \cite{ARS1}]\label{thm16.1}
 Let $V$ be a closed invariant subspace of $\bbh = A^2(\bbc, D^1(0, 1), \dmu)$, where $\dmu$ is real $2$-dimensional area measure $\dx \dy $.. Then $V$ is the smallest closed invariant subspace generated by $M$\!, where $M = V \orthcomp z V$\!.
\exitthmcall{theorem}

The ARS theorem is different, and arguably better, than ours. The apparent paradox it raises involves the dimension of $M$.  Let

\begin{equation*}
   V = W_0 \plusperp W_1 \plusperp W_2 \plusperp \cdots  
\end{equation*}
\noindent be the valuation homogeneous decomposition of $V$\!. Recall that by the third projection lemma, lemma \ref{lem12.6}, the projection operators

\begin{equation*}
    L_{ W_m }  \colon W_m \mapsto P_m^\bbh ( W_m ) \subseteq H_m, \,\,m = 0, 1, 2, \ldots
\end{equation*}
\noindent are invertible.  Now in $1$ complex variable, each subspace $H_m = \LH ( z^m ), m = 0, 1, 2, \ldots $ has dimension $1$ (or $0$).  In other words, $ \bbh $ is singly generated as defined in definition \ref{def18.4} below. Thus each subspace

\begin{equation*}
      W_m = V_m \orthcomp \! \! V_{m + }, \, \,  m =  1, 2, 3, \ldots
    \end{equation*}
\noindent has dimension $0$ (if it equals $\{ 0 \}$) or $1$. But, as Aleman, Richter, and Sunberg showed, $M$ can have any dimension from $1$ up to and including $\infty$.

\

The resolution of this paradox is contained in a general observation, which we shall state momentarily as a proposition. Recall that $V_k = \{ h \in V \colon \ord(h) \geq k \}$.  Let $k$ be the first index such that $W_k \neq \{ 0 \}$, then $V = V_{ k }$. Thus $M = V_{ k } \orthcomp z V_{k}$. Our general observation is,

\thmcall{proposition}\label{prop16.2}
    Let $ \bbh $ be an analytic Hilbert module over the analytic algebra $ \bbp [ \bbc, z ] $. Let $ k $ be the first index index such that $ W_k \neq \{ 0 \} $. Suppose $ M $ has dimension $ \geq 2 $. Then $ z V_k $ is a proper subset of $ V_{ k+ 1} $, so $ M $ is a proper superset of $ W_{ k + 1 } $.
\exitthmcall{proposition}

\begin{proof}
    If $ h \in V_k $, then $ \ord( h ) \geq k $, so $ \ord ( z h ) \geq k+ 1 $. Thus $ z V_k \subseteq V_{k + 1} $, so $ M = V_k \orthcomp z V_k  \supseteq  V_k \orthcomp V_{ k + 1 } = W_k $. Now if $ M $ has dimension $ \geq 2 $, then $ M \supsetneq W_k $, so $ z V \subsetneq V_{ k + 1 } $.
\end{proof}

The difficulty with this proof is that it appears to be circular, in that it appears to use the paradox to resolve the paradox. Here is more direct proof for $A^2(\bbc, D^1(0, 1), \dmu)$, simply using the power series expansions of analytic functions in $A^2(\bbc, D^1(0, 1)), \dmu)$:

\begin{proof}
    As above, $ z V_k \subseteq V_{k + 1} $. Now suppose $ z V_k $ actually equals $ V_{ k + 1 } $ and
\begin{equation*}
     M =  V_k \orthcomp z V_{k+1} =  V_k \orthcomp V_{k+1}
\end{equation*}
so $M \perp V_{k+1}$.

   Let $ f $ and $ g $ be two linearly independent elements of $M$. Note that $ f $ and $ g \in V_k $. Let $ f ( z ) = a_k z^k + a_{k + 1} z^{k + 1}  + a_{k + 2} z^{k + 2} +\cdots $, and $ g ( z ) = b_k z^k + b_{k + 1} z^{k + 1} + b_{k + 2} z^{k + 2} + \cdots $. There are two cases:
    
    If both of $ a_k $ and $ b_k \neq 0 $, then, by multiplying $ g $ by a constant if necessary, without loss of generality, $ a_k = b_k $. Let $ h = f - g $, and let $ h ( z ) = c_k z^k + c_{k + 1} z^{k =1} + c_{k + 2} z^{k + 2} + \cdots $.  Then $ h \neq 0 $ because $ f $ and $ g $ are linearly independent, $ h \in M $ because $ M $ is a linear subspace, and $ c_k = a_k - b_k = 0 $, so $ h \in V_{ k + 1 } $.  But this is impossible because $M \perp V_{k+1}$.
   
    On the other hand, suppose $ a_k =  0$ or $b_k = 0$.   Set $ h = f $ if $a_k = 0$, or $h = g$ if $b_k = 0$.  As above, $ h \neq 0 $, $h \in M$, and $ c_k = 0 $, so $ h  \in V_{ k + 1 } $, which is, as we have just seen, impossible.
\end{proof}

Hedenmalm, Richter, and Seip have given an explicit construction of cases, involving the angle between subspaces, where $M = V \orthcomp zV$ has finite or infinite dimension. Their construction is valid, not only for $A^2(\bbc, D^1),dx \dy$, but also for weighted Bergman spaces in the unit disk \cite{HRS1}.

Walter Rudin has given an example of a closed invariant subspace $V$ of $H^2(\bbc, D^2, \\ \dx \dy)$ (with $D^2 = D^2(0, 1) = $ the unit polydisk in $\bbc^2$ and $R = P_2[\bbc, z_1, z_2]$) which is not finitely generated (\cite{Ru2}, Corollary to his Theorem 4.4.2). While this might not be too surprising, it does show that being infinitely generated is as much the rule as the exception for $V  \orthcomp R_1 V$ for analytic Hilbert modules in several complex variables.

\section{$ R_1 $ Invariant Subspaces and $ R_1 $ Inner Functions.}\label{sec17}

This section will be devoted to the study of $ R_1 $ inner functions and their relation to closed $ R_1 $ invariant subspaces.  We shall show that every non-zero closed $ R_1 $ invariant subspace of a valuation Hilbert module is generated by a collection of $ R_1 $ inner functions.

Throughout this section, $ R $ will be the valuation algebra $ \bbp_n = \bbp [ \bbc, z_1, \ldots, z_n] $.  $ \bbh $ will be one of the valuation Hilbert modules
$ H^2( \bbc, D^n ( p, r ) ) $, $ H^2 ( \bbc, B^n ( 0, 1 ) ) $,\\ or  $ A^2 ( \bbc, \Omega, k\dmu ) $, all over $ R $, where $ \Omega $ is a domain in $ \bbc^n $ and $\dmu$ is real $2 n$-dimensional volume measure on $\Omega$, ${\dx}_1 {\dy}_1 \cdots {\dx}_n {\dy}_n$.

And throughout this section, unless the contrary is explicitly stated, $V$ will be a closed non-zero subspace of $ \bbh $.  $ V = W_0 \plusperp W_1 \plusperp W_2 \plusperp \dots $ will, as usual, be the valuation homogeneous decomposition of $ V $.  And $  \bbh = H_0 \plusperp H_1 \plusperp H_2 \plusperp \, \cdots $ will be the \emph{valuation homogeneous decomposition} of $ \bbh $ itself.

Recall from definition \ref{def10.1} that a closed non-zero subspace W of $ \bbh $ is \emph{$ R_1 $ inner} if $ R_1 \cdot W \perp W $. Recall also that the definition of a valuation homogeneous decomposition which is $ R_1 $ \emph{inner} is more complicated because it involves the essential ``look ahead" condition:  If, for $k$ and $m \in \bbz_+$, $h \in W_k$, $g \in V$, $r \in R_1$, and $\ord(r h - g) > m$, then $ r h - g \perp W_m $.

\thmcall{definition} \label{def17.1}
    Similarly, a non-zero function $ h $ is \emph{ $ R_1 $ inner} if the closed one dimensional subspace $ \LH (h ) = \{ \lambda \cdot h \colon \lambda \in \bbc \} $ is $ R_1 $ inner.
\exitthmcall{definition}

Thus our abstract Beurling's theorem could be rephrased to say that every non-zero closed invariant subspace is generated by a collection of $ R_1 $ inner functions (with the addition, of course, of the look ahead condition where needed).

\section{The Characterization of $ R_1 $ Inner Functions}\label{sec18}

In this section, we shall give examples showing that the most obvious characterizations of $ R_1 $ inner functions are false.

Throughout this section, $ R $ will be the algebra 
 
 \begin{equation}\label{eqn18.1}
    \bbp_n = \bbp_n[ \bbc, z ]  = \bbp [ \bbc, z_1 \ldots z_n ]
\end{equation}
of complex-valued polynomials in the complex variables $ z = z_1, \ldots z_n $, and $ \bbh $ will be one of the the simplest Hilbert modules over $ R $, namely the norm closure of $ \bbp_n $. Of course, to form $\bbh$, we must define an inner product and norm on $ \bbp_n $.

Another useful way to view $ \bbh $ is

\begin{equation}\label{eqn18.2}
\begin{split}
   \bbh =  &\text{ the scalar valued Hilbert space generated by the set of monomials} \\
               &\,\,\{ z^k = z_1^{ k_1} \ldots, z_n^{k_n} \colon k \in \bbz_{+}^n \}. \\
\end{split}
\end{equation}

\thmcall{remark}\label{rem18.3}
To begin, we shall construct the valuation Hilbert module $ \bbh $ in a general way:

First, we shall construct an inner product $ ( h_1, h_2) \mapsto \inprod{ h_1 }{ h_2 } $ and a norm $ \norm{ h }^2 = \inprod{ h }{ h } $ on $ \bbp_n $.

Second, we shall construct the inner product so the monomials in the set of monomials in equation \ref{eqn18.2} are orthogonal.

Third, we shall take the norm completion of $ \bbp_n $ to form the complex Hilbert space $ \bbh $.  Then the set of normalized monomials from equation \ref{eqn18.2} will form a complete orthonormal system for $ \bbh $.

Fourth, and finally, we shall let the $ \ord $ functions on both $ R $ and $ \bbh $ be standard order functions.  Then $ ( R, \ord ) $ and $ ( \bbh, \ord ) $ will be, respectively, an analytic valuation algebra and a valuation Hilbert module over $ R $.
\exitthmcall{remark}

While this might seem to be overly restrictive, most of the valuation Hilbert modules we have studied so far are closed submodules, that is closed $ R_1 $ invariant subspaces of Hilbert modules over $ R = \bbp_n $ of this form.  These include all but one of the $ H^2 $ modules studied in section \ref{sec14}, and weighted Bergman spaces on circularly symmetric complex domains with circularly symmetric weight functions, which are special cases of those studied in section \ref{sec15}.

\

We need a name for such valuation Hilbert modules:

\thmcall{definition} \label{def18.4}
A valuation Hilbert module $ \bbh $ is \emph{$ \bbp_n $ generated} if it is the norm closure of the linear subspace $ \bbp_n $, and if the set of monomials \ref{eqn18.2} form a complete orthogonal system for $ \bbh $.

We shall also say such valuation Hilbert modules are \emph{polynomially generated} when we do not need to keep track on the dimension $ n $.

If $ n = 1 $, so $ \bbh $ is $ \bbp_1  $ generated, we shall say that \emph{$ \bbh $ is singly generated}.
\exitthmcall{definition}

\thmcall{remark} \label{rem18.5}
In what follows, we shall freely mix functional analysis notation $ h $ and hard analysis notation $ h ( z ) $ when referring to a function $ h $.  It will be clear from context whether $ h ( z ) $ refers to the function $ h $ or the value of $ h $ at the point $ z $. 
\exitthmcall{remark}

First, we shall consider the valuation Hilbert module $ \bbh = H^2 ( \bbc, D^n ( 0, 1 ) $ over the valuation algebra $ R = \bbp_n $.  $ T^n $ will be the distinguished boundary of $ D^n(0, 1) $.  Of course, $ T^n $ is the $n$-dimensional unit torus, $ z = \{ (z_1, \ldots, z_n) \} $ with $ \abs{z_j} = 1 $ for $j = 1, \ldots, n$.

\thmcall{remark}\label{rem18.6}
We need to remind the reader of a few well known facts about $ \bbh $ and $ T^n $.  These can be found in chapter 3 of \cite{Ru2}.

\enumcall{enumerate}
  \item The non-tangential boundary values $ f^* $ of $ f \in \bbh $ exist a.e.\@ $ \dtheta $ on $ T^n $, where, as usual,
    \begin{equation*}
       \begin{split}\label{item18.6.1}
          \theta &= \theta_1, \ldots, \theta_n. \\
          \dtheta &= \dtheta_1 \cdots \dtheta_n. \\
       \end{split}
    \end{equation*}
  \item The measure $ \dtheta $ has total mass $ ( 2 \pi )^n $ on $ T^n $. \label{item18.6.2}
  \item The set of monomials in equation \ref{eqn18.2}, $ \{ z^k \colon  k \in \bbz_+^n \} $,  form a complete orthonormal system for $ \bbh $. \label{item18.6.3}
\exitenumcall{enumerate}
\exitthmcall{remark}

Let $ n > 1 $ and $ f \neq 0 $ be an $ R_1 $ inner function in $ H^2 ( D^n ( 0, 1 ) ) $.  By analogy with the one complex variable case, one might think that $ \abs{ f^* } $ must be constant a.e.\@ on $ T^n $.  But the simplest examples show that this is not so:

\thmcall{example}\label{example18.7}
   Let $ n = 2 $, and consider $ V = \bbh = H_0 \plusperp H_1 \plusperp H_2 \plusperp \cdots $.  Because this expansion is the valuation homogeneous decomposition of $ \bbh $ and $ \bbh $ is itself $ R_1$ invariant, each $ H_k $ is $ R_1 $ inner.  Let $ f(z ) = f ( z_1, z_2 ) = z_1 \cdot z_2 +  z_2^2 $.  Because $ z $ is continuous on $ \bbc^n $, we may omit the $ * $ superscripts.

   $ f (1, 1 ) = 1 \cdot 1 + 1^2 = 2 $.  Because $ f $ is continuous, $ \abs{ f } > 1 $ on an open subset of $ T^n $ containing $ ( 1, 1 ) $.  Thus $ \abs{ f } > 1 $ on a subset of $ T^n $ of positive measure containing the point $ ( 1, 1 ) \in T^n $.

   On the other hand, $ f ( 1, -1 ) = 1 \cdot ( -1 ) +  ( -1 )^2 = ( -1 ) + 1 = 0 $.  As above, $ \abs{ f } < 1/2 $ on a subset of $ T^n $ of positive measure containing the point $ ( 1, -1 ) \in T^n $.

   Hence $ \abs{ f } $ is not constant a.e.\@ on $ T^n $.
\exitthmcall{example}

Next, we shall consider the unit ball $ B^n = B^n ( 0, 1 ) $ in $ \bbc^n $, and the valuation Hilbert module $ \bbh = H^2 ( \bbc^n, B^n ) $ over the valuation algebra $ R = \bbp_n $.  By contrast with the polydisk, the distinguished boundary of $ B^n $ coincides with its topological boundary $ \bdry B^n $, and equals the real $ 2n -1 $ dimensional sphere $ S^{ 2n - 1 } = S^{ 2n - 1 } ( 0, 1 ) $ of radius $ 1 $ centered at $ 0 $.

\thmcall{remark}\label{rem18.8}
We need to remind the reader of a few well known facts about $ \bbh $.  These can be found in \cite{Ru4}, especially in chapter 5.

\enumcall{enumerate}
  \item The non-tangential boundary values $ f^* $ of $ f \in \bbh $ exist a.e.\@ $ d\sigma $ on $ S^{ 2n - 1 } $, where $ d\sigma $ is normalized real $ 2n - 1 $ volume measure on $ S^{ 2n -1 } $.\label{item18.8.1}
  \item In multi-variable notation, for $ k \in \bbz_+^n $, $ 0 < \norm{ z^k } < \infty $.  Of course, $ \norm{ z^k }^2 = \inprod{ z^k }{ z^ k } $. \label{item18.8.2}
  \item The set of monomials in equation \ref{eqn18.2}, $ \{ z^k \colon  k \in \bbz_+^n \} $,  form a complete orthonormal system for $ \bbh $. \label{item18.8.3}
\exitenumcall{enumerate}
\exitthmcall{remark}

\thmcall{example}\label{example18.9}
   As in the case of the polydisk, the $ R_1 $ inner function $ f(z ) = f ( z_1, z_2 ) = z_1 \cdot z_2 +  z_2^2 $ in $ \bbh $ is not constant on the distinguished boundary $ S^{ 2n - 1 } $.  At the North pole, $ f ( 1, 0 ) = 0 $.  At the East pole, $ f ( 0, 1 ) = 1 $.
\exitthmcall{example}

We turn now to the unweighted Bergman space $ \bbh = A^2 ( \bbc, B^n ( 0, 1 ), \dmu ) $, where $ \dmu $ is normalized real $ 2n $-dimensional volume measure on $ B^n = B^n ( 0, 1 ) $, and $ R = \bbp_n $.  A transformation to polar coordinates shows that the measure $ \dmu $ is, in multi-index notation, 
\begin{equation} \label{eqn18.10}
  \dmu = \frac{ n! }{ \pi^n } \cdot dx dy
\end{equation}
(where $ \dx \dy $ is real $ 2n $-dimensional volume measure) and the Bergman kernel for $ B^n $ is,
\begin{equation} \label{eqn18.11}
  K_n ( z, w ) =  \sum \,\frac{ z^k }{ \norm{ z^k } } \cdot \frac{ \overline{ w }^k }{ \norm{ \overline{ w }^k } } = \frac{ n! }{ \pi^n } \cdot \frac{ 1 }{ (1 - z \cdot \overline{ w} )^{ n + 1 } }
\end{equation}
(Cf.\@ \cite{Kr1}, theorem 1.4.21.)  Thus, in multi-index notation, the normalized monomials $ z^k /\norm{ z^k }, \,k \in \bbz+^n $ form a complete orthonormal system for $ \bbh $.  Thus we have

\thmcall{example}\label{def18.12}
As in the cases of the $ H^2 $ of the polydisk and unit ball, some $ R_1 $ inner functions in $ \bbh $ do not have modulus $ = 1 $ ae.\@ on the boundary $ S^{ 2n - 1 } $.

In fact, even in dimension $ 1 $ for the unit disk, the situation can be bad.  Some analytic functions in $ \bbh = A^2 ( \bbc, D^1 ( 0, 1 ), \dmu ) $ \emph{fail to have} radial boundary values \emph{everywhere} on the unit circle \cite{Hed1}.  It is unclear whether some $ R_1 $ inner functions also lack non-tangential or radial boundary values a.e.\@ on the unit circle.
\exitthmcall{example}

\section{Standard Inner Functions.} \label{sec19}

\thmcall{definition}\label{def19.1}
Recall that a complex-valued function $ I $ in $ H^2 ( \bbc, D^1 ( 0, 1 ) ) $ of the unit disk $ D^1 ( 0, 1 ) $ is normally called an \emph{inner function} if it has modulus $ = 1 $ a.e.\@ on the unit circle, $ T^1 $.  To distinguish these from our $ R_1 $ inner functions, we shall call such functions \emph{standard inner functions} (not to be confused with \emph{standard order functions}).

Note that the usual definition of an inner function may be extended to complex-valued functions in $ H^2 ( \bbc, B^n ( 0, 1 ), \dmu ) $ of the $ n $-dimensional unit ball $ B^n ( 0, 1) $, except that the boundary in this case is the unit real $ 2n - 1 $ dimensional sphere $ S^{ 2n - 1 } $.  Again, we shall call such functions \emph{standard inner functions}.
\exitthmcall{definition}

For many years, it was unclear whether non-constant, standard inner functions in $ \bbh = H^2 $ of the $ n $-dimensional unit ball existed..  The deep conjecture that such standard inner functions did not exist came to be known as the \emph{Rudin Inner Functions Conjecture} because, in section 19.1 of \cite{Ru4}, Walter Rudin gave many examples of the pathological properties that standard inner functions in $ \bbh $ would have to possess.  It thus came as a great surprise that such functions did exist.  (Cf.\@ the discussions in \cite{Ru5} and chapter 9 of \cite{Kr1}.)  Credit for this discovery is due independently, and by somewhat different methods, to four authors, Aleksandrov, Hakim and Sibony, and Low.  (Cf.\@ the discursive discussion in \cite{Ru5}.)

The point is that standard inner functions in many cases are far more complicated than our $ R_1 $ inner functions. 

The original extension of Beurling's theorem to $A^2$ spaces on the unit disk by Aleman, Richter, and Sunberg involved biharmonic functions and was notoriously difficult (\cite{ARS1}, \cite{Hed1}). Since then, Ball and Bolotnikov have given a different formulation of Beurling's theorem for these spaces for one complex variable involving standard inner functions (\cite{BaBo1}, \cite{BaBo2}). It would be of great interest to extend Ball and Bolotnikov's formulation to several complex variables.

However, our formulation of Beurling's theorem for $A^2$ spaces involving valuation Hilbert modules, even in one complex variable, has the signal advantage that the $R_1$ inner functions involved are \emph{complex-valued}. This is particularly useful when considering questions surrounding the \emph{invariant subspace conjecture} (cf.\@ \cite{HRS1} and \cite{HKZ1},  page 187)). The reader might recall that de Branges' and Rovniak's 1964 proof of the invariant subspace conjecture failed because their method for factoring \emph{operator-valued} inner functions was flawed \cite{dBr1}, \cite{dBr2}.  Methods, like those in our paper involving \emph{complex-valued} functions, may allow one to bypass the deep problem of factoring \emph{operator-valued} inner functions.

\section{$ A^p_\alpha $ and $ R_1 $ Inner Functions.} \label{sec20}

Hedenmalm, Korenblum, and Zhu present an extensive theory of \emph{$ A^p_\alpha $ inner functions} in weighted Bergman spaces $ A^p_\alpha $ in the unit disk.  Here $ p $ refers to $ L^p $ and $ \alpha  $ refers to the weight $ ( 1 - \abs{ z }^2 )^\alpha $ (cf.\@ \cite{HKZ1}, p.~2 and chapter 3).  In our notation, $ A^p_\alpha $ in the unit disk for $ p = 2 $ is $ A^2( \bbc, D^1(0, 1), k \dmu) $, where $ k $ is the above weight.

As Hedenmalm, Korenblum, and Zhu point out, an $ A^p_\alpha $ inner function is the direct analog of a standard inner function in $ H^p $ of the unit disk because, to paraphrase them slightly, and in our notation for $ p = 2 $,

\thmcall{definition} \label{def20.1}
A function $ f $ in $ A^2_\alpha $ is an $ A^2_\alpha $ inner function if and only if
   \begin{equation} \label{eqn20.1.1}
      \int_{ D^1 ( 0, 1 ) } \abs{ f ( z ) }^2 \,r ( z ) ( 1 - \abs{ z }^2 )^\alpha \dx\dy / \pi = r (0 ).
   \end{equation}
for every polynomial $ r \in \bbp_1 = \bbp [ \bbc, z ] $.  For $ H^2 $ of the unit disk, this condition (with $ \dtheta / 2 \pi $ in place of $ \dx\dy / \pi $) is equivalent to $ \abs{ f } $ being constant ae.\@ on the unit circle. (Cf.\@ \cite{HKZ1}, definition 3.1.)
\exitthmcall{definition}

\thmcall{proposition} \label{prop20.2}
$ A^2_\alpha $ inner functions and our $ R_1 $ inner functions in $ A^2_\alpha $ are the same for $ R_1 = \bbp_{1, 0} $, the algebra of polynomials $ r $ in $ 1 $ variable with $ r ( 0 ) = 0 $.
\exitthmcall{proposition}

\begin{proof}
Let $ h $ be an $ A^2_\alpha $ inner function.  From definition \ref{def20.1}, we see that $ z^m \inprod{ h }{h} = 0 $ for $ m = 1, 2, 3 \ldots $.  Thus $ \inprod{ z^m h }{ h } =  0 $ for $ m = 1, 2, 3, \ldots $, so $ \inprod{ r h }{ h } = 0 $ for all $ r \in R_1 $.  Consequently, $ h $ is an $ R_1 $ inner function.  Conversely, suppose that $ h $ is an $ R_1 $ inner function.  The same argument shows that $ h $ is an $ A^2_\alpha $ inner function.
\end{proof}

\section{Open Questions.}\label{sec21}

\thmcall{question}\label{quest22.1}
     In several complex variables, how can we characterize the closed invariant subspaces of $H^p(E, \Omega)$ and $A^p(E, \Omega,$ \,$k\dmu)$ for $1 \leq p <  \infty$, and of $H^\infty(\Omega)$ under the $\beta$ topology, that is under Buck's $\beta$ topology of bounded compact convergence?  Here, $E$ is a complex Hilbert space.  See Buck \cite{Bu1} for his $\beta$ topology.  In the one complex variable case for $E = \bbh$ and $\Omega =$ the unit disk, the classical Beurling's theorem generalizes to these $H^p$ spaces.  (See \cite{Hel1}, p.\@ 25.)  For $E$ = a complex Hilbert space and $p = 2$, see \cite{Hoff1}.  See \cite{RS1} and \cite{RS2} for complex-valued $H^\infty$ on the unit disk under the $\beta$ topology, and \cite{N1} and \cite{N2} for a generalization to these complex-valued $H^p$ spaces on open Riemann surfaces.
\exitthmcall{question}

\thmcall{question}\label{quest21.2}
    For what dimensions $n$, for what domains or even complex analytic manifolds $ \Omega $, and for what weighted $ A^2 $ spaces $ A^2( \bbc, \Omega, k\dmu) $, can we deduce the ARS characterization of closed invariant subspaces, $ V = \cl R_1 \cdot M $, where $ R = \bbp_n [ \bbc, z_1, \ldots z_n ] $ or $ H^\infty(\Omega) $, and $ M = V \orthcomp R_1 \cdot V $, from our abstract Beurling's theorem?
\exitthmcall{question}

Of course, we can ask question \ref{quest21.2} even more generally:

\thmcall{question}\label{quest21.3}
    For what analytic algebras $R$ and analytic Hilbert modules $\bbh$, including vector valued ones, can we deduce the ARS characterization of closed invariant subspaces?
\exitthmcall{question}

\thmcall{question}\label{quest21.4}
    Halmos's paper also characterized the closed invariant subspaces of $L^2$ of the unit circle. See the treatments in the books by Helson \cite{Hel1} and Hoffman \cite{Hoff1}. How can we generalize our abstract Beurling's theorem to the case of $L^2$ spaces in several complex variables?  The same question for $L^p$ spaces in several complex variables.
\exitthmcall{question}

\thmcall{question}\label{quest21.5}
    The original extension of Beurling's theorem to $A^2$ spaces on the unit disk by Aleman, Richter, and Sunberg involved biharmonic functions (\cite{ARS1}, \cite{Hed1}). Since then, Ball and Bolotnikov have given a different formulation of Beurling's theorem for these spaces for one complex variable involving standard inner functions (\cite{BaBo1}, \cite{BaBo2}). It would be of great interest to extend Ball and Bolotnikov's formulation to several complex variables.
\exitthmcall{question}

\thmcall{question}\label{quest21.6}
    What is the connection between operator valued standard inner functions and our $ R_1 $ inner functions?  By contrast with operator-valued standard inner functions, the $ R_1 $ inner functions we considered in sections \ref{sec17} and \ref{sec18} have the signal advantage that they are complex-valued.

    As we observed in the discussion at the end of section \ref{sec19}, this might be particularly useful when considering questions surrounding the invariant subspace conjecture. See Hedenmalm, Richter, and Seip's paper, \cite{HRS1}, and Hedenmalm, Korenblum, and Zhu's book, \cite{HKZ1}, page 187, for a discussion of their reduction of the invariant subspace conjecture to a question about invariant subspaces of $A^2$ of the unit disk.
    
    The reader might also recall that de Branges' and Rovniak's 1964 proof of the invariant subspace conjecture failed because their method for factoring operator-valued inner functions was flawed \cite{dBr1}, \cite{dBr2}.  Methods like those in our paper involving complex-valued $R_1$ inner functions in valuation Hilbert modules might allow one to bypass the deep problem of factoring operator-valued inner functions.
\exitthmcall{question}


\begin{thebibliography}{99}

\bibitem{ACD1}
O.~P. Agrawal, D.N. Clark, and R.~G Douglas, \emph{Invariant subspaces in the
  polydisk}, Pacific J. Math. \textbf{121} (1986), no.~1, 1--11. MR 815027  Zbl 0609.47012
  
\bibitem{AC1}
P.~R. Ahern and D.~N. Clark, \emph{Invariant subspaces and analytic
  continuation in several variables}, J. Math. Mech. \textbf{19} (1969/1970),
  963--969. MR 0261340 Zbl 0211.10203

\bibitem{AB1}
S.~Axler and P.~Bourdon, \emph{Finite-codimensional invariant subspaces of
  {B}ergman spaces}, Trans. Amer. Math. Soc. \textbf{306} (1988), no.~2,
  805--817. MR 933319 Zbl 0658.47011
  
\bibitem{ARS1}
A.~Aleman, S.~Richter, and C.~Sundberg, \emph{Beurling's theorem for the
  {B}ergman space}, Acta Math. 177 (1996), no.~2, 275--310.
  MR 1440934  Zbl 0886.30026

\bibitem{BaBo1}
J.~A.~Ball and V.~ Bolotnikov, \emph{A {B}eurling type theorem in
  weighted {B}ergman spaces}, C. R. Math. Acad. Sci. Paris 351 (2013),
  no.~11-12, 433--436. MR 3090124 Zbl 1277.47012

\bibitem{BaBo2}
\bysame
\emph{Weighted {B}ergman spaces:~shift-invariant subspaces and input/state/output linear systems}, Integral Equations and Operator Theory 76 (213), no.~3, 301-356. MR 3065298 Zbl 

\bibitem{Beu1}
A.~Beurling, \emph{On two problems concerning linear transformations in
  {H}ilbert space}, Acta Math. \textbf{81} (1948), 17. MR 0027954 Zbl 0033.37701
  
\bibitem{dBr1}
L.~de~Branges and J.~Rovniak, \emph{The Existence of Invariant Subspaces}, Bull. Amer. Math. Soc. 70 (1964), 721-721.  Not reviewed. 

\bibitem{dBr2}
\bysame, \emph{Correction to ``The existence of invariant subspaces"}, Bull. Amer. Math. Soc. 71 (1965), 396.  MR  0169024  Zbl 0171.34902

\bibitem{Bu1}
R.~Buck, \emph{Algebraic Properties of Classes of Analytic Functions}, Seminars on Analytic Functions, vol. 2, Princeton University Press (1957), 175--188. MR not available. Zbl 0196.43603
  
\bibitem{DPSY1}
R.~G. Douglas, V.~I. Paulsen, C.~H. Sah, and K.~Yan, \emph{Algebraic reduction
  and rigidity for {H}ilbert modules}, Amer. J. Math. \textbf{117} (1995),
  no.~1, 75--92. MR 1314458  Zbl 0833.46040
 
 \bibitem{Gr1}
A.~ Grigor$\prime$yan, \emph{Analytic and geometric background of
  recurrence and non-explosion of the {B}rownian motion on {R}iemannian
  manifolds}, Bull. Amer. Math. Soc. (N.S.) 36 (1999), no.~2,
  135--249. MR 1659871  Zbl 0927.58019
  
\bibitem{Guo1}
K.~Guo, \emph{Characteristic spaces and rigidity for analytic {H}ilbert
  modules}, J. Funct. Anal. \textbf{163} (1999), no.~1, 133--151. MR 1682835 Zbl 0937.46047

\bibitem{Hal1}
P.~R. Halmos, \emph{Shifts on {H}ilbert spaces}, J. Reine Angew. Math.
  \textbf{208} (1961), 102--112. MR 0152896  Zbl 0107.09802

\bibitem{Has1}
M.~Hasumi, \emph{Invariant subspaces on open {R}iemann surfaces}, Ann.
  Inst. Fourier (Grenoble) 24, no.~4, (1974) 241--286. MR 0364647 Zbl 0287.46066

\bibitem{Has2}
\bysame, \emph{Invariant subspaces on open {R}iemann surfaces. II}, Ann.
  Inst. Fourier (Grenoble) 26, no.~4,  (1976) 273--299. MR 0407283 Zbl 0287.46066

\bibitem{Hed1}
H.~Hedenmalm, \emph{Recent progress in the function theory of the {Bergman} space}, in \emph{Holomorphic spaces}, MSRI Publications, 33 (1998) 35--50, Cambridge University Press. MR 163064 Zbl 1126.46303.  Also at http://library.msri.org/books/Book33/files/heden.pdf, accessed Thursday, July 29, 2022.

\bibitem{HRS1}
H.~Hedenmalm, S.~Richter, and~K.~Seip,~\emph{Interpolating sequences and invariant subspaces of given index in the {B}ergman spaces}, J. Reine Angew. Math., 477 (1996) 13--30.  MR 1405310 Zbl 0895.46023.  Also at https://eudml.org/doc/153837$\ast$, accessed Friday, December 28, 2020.

\bibitem{HKZ1}
H.~Hedenmalm, B~Korenblum, and K.~Zhu,~\emph{Theory of Bergman Spaces}, Graduate Texts in Mathematics, No.~199, Springer Verlag, Berlin-NewYork, 2000.  MR1758653 Zbl 0955.32003 

\bibitem{H1}
M.~Heins, \emph{Selected topics in the classical theory of functions of a
  complex variable}, Athena Series: Selected Topics in Mathematics, Holt,
  Rinehart and Winston, New York, 1962. MR 0162913 Zbl 1226.30001
 
 \bibitem{H2}
\bysame, \emph{Hardy classes on {R}iemann surfaces}, Lecture Notes in
  Mathematics, No.~98, Springer-Verlag, Berlin-New York, 1969. MR 0247069 Zbl 

\bibitem{Hel1}
H.~Helson, \emph{Lectures on invariant subspaces}, Academic Press, New
  York-London, 1964. MR 0171178  Zbl  0119.11303
  
\bibitem{Hoff1}
K.~Hoffman, \emph{Banach spaces of analytic functions}, Prentice-Hall Series in
  Modern Analysis, Prentice-Hall, Inc., Englewood Cliffs, N. J., 1962.
  MR 0133008 Zbl 0117.34001

\bibitem{Kr1}
S~Krantz, \emph{Function theory of several complex variables}, second
  ed., AMS Chelsea Publishing, 2001. MR 1846625 Zbl 1087.32001
 
\bibitem{MASS1}
 A.~Maji, A.~ Mundayadan, J.~Sarkar, and T.~R.~Sankar, \emph{Characterization of invariant subspaces in the polydisc}, Preprint: arXiv:1710.09853v2,  (2017)

\bibitem{MASS2}
 A.~Maji, A.~ Mundayadan, J.~Sarkar, and T.~R.~Sankar, \emph{Characterization of invariant subspaces in the polydisc}, J. Operator Theory 82 (2019), no. 2, 445-468. MR 4015959 Zbl 1438.47020

\bibitem{N1}
C.~W.~Neville, \emph{Invariant subspaces of {H}ardy classes on infinitely
  connected plane domains}, Bull. Amer. Math. Soc. 78 (1972),
  857--860. MR 0301206   Zbl 0266.46040

\bibitem{N2}
\bysame, \emph{Invariant subspaces of {H}ardy classes on infinitely connected
  open surfaces}, Mem. Amer. Math. Soc. 2 (1975), issue 1, 160,
  viii+151. MR 0586558   Zbl 0314.46052
  
 \bibitem{N3}
\bysame, \emph{An abstract {B}eurling's theorem for several complex variables I}, Preprint 2021, arXiv:2108.12272, 23 pages.

\bibitem{N4}
\bysame, \emph{An abstract {B}eurling's theorem for several complex variables II}, Preprint 2021, arXiv:2109.00695, 17 pages.

\bibitem{N5}
\bysame, \emph{{B}eurling's theorem for valuation  {H}ilbert modules and several complex variables}, Preprint 2022, arXiv:2109.00695, 25 pages.

\bibitem{Pu1}
M.~Putinar, \emph{On invariant subspaces of several variable {B}ergman spaces},
  Pacific J. Math. \textbf{147} (1991), no.~2, 355--364. MR 084714  Zbl 0692.47008

\bibitem{RS1}
L.~A.~Rubel and A.~L.~Shields, \emph{Weak topologies on the bounded holomorphic functions},  Bull. Amer. Math. Soc. 71 (1965), 349--352.  MR 172138  Zbl  0173.41601

\bibitem{RS2}
L.~A.~Rubel and A.~L.~Shields, \emph{The space of bounded holomorphic functions on a region},  Ann. Inst. Fourier (Grenoble) 16 (1966), no. 1, 235--277.  MR 1533135  Zbl  0152.13202

\bibitem{Ru1}
W.~Rudin, \emph{Analytic functions of class {$H_p$}}, Trans. Amer. Math.
 Soc. 78 (1955), 46--66. MR 0067993  Zbl 0064.31203

\bibitem{Ru2}
\bysame, \emph{Function theory in polydiscs}, W. A. Benjamin, Inc., New
  York-Amsterdam, 1969. MR 0255841  Zbl 0177.34101

\bibitem{Ru3}
\bysame, \emph{Functional analysis}, McGraw-Hill Series in Higher Mathematics, McGraw-Hill, Inc.\@, New York, 1973.

\bibitem{Ru4}
\bysame, \emph{Function theory in the unit ball of {${\bf C}^{n}$}},
  Grundlehren der Mathematischen Wissenschaften [Fundamental Principles of
  Mathematical Science], vol. 241, Springer-Verlag, New York-Berlin, 1980.
  MR 601594  Zbl 0495.32001
 
\bibitem{Ru5}
\bysame, \emph{New constructions of functions holomorphic in the unit ball of {${\bf C}^{n}$}}, CBMS Regional Conference Series in Mathematics, vol. 63,  American Mathematical Society, Providence, RI, 1986. xvi+78 pp. MR 840468 Zbl 1187.32001
 
\bibitem{Sh1}
S.~Shimorin, \emph{Wold-type decompositions and wandering subspaces for
  operators close to isometries}, J. Reine Angew. Math. 531 (2001),
  147--189. MR 1810120 Zbl 0974.47014

\bibitem{Sh2}
\bysame, \emph{On {B}eurling-type theorems in weighted {$l^2$} and
  {B}ergman spaces}, Proc. Amer. Math. Soc. 131 (2003), no.~6,
  1777--1787. MR 1955265 Zbl 1046.47024

\bibitem{SY1}
M.~Seto and R.~Yangi, \emph{Inner sequence based invariant subspaces in
  {$H^2(D^2)$}}, Proc. Amer. Math. Soc. \textbf{135} (2007), no.~8, 2519--2526.
  MR 2302572  Zbl 1137.47006
 
\bibitem{Wiki1}
Wikipedia Contributors, \emph{Semi-continuity}, Wikipedia, the free
  encyclopedia (2021), 1, Accessed February 20, 2021. Not indexed by MR or Zbl.
  
\end{thebibliography}
\end{document}